                      \def\version{October 08, 2009}                           % 
\documentclass[reqno,12pt]{amsart} 
\usepackage{amsmath} 
\usepackage{amssymb} 
\newenvironment{proofsect}[1] 
{\vskip0.1cm\noindent{\bf #1.}\hskip0.5cm}

\def\emptyset{\varnothing} 
\def\ti{\to\infty}

\def\d{\delta}

\def\L{\Lambda}

\newfam\Bbbfam 
\font\tenBbb=msbm10 
\font\sevenBbb=msbm7 
\font\fiveBbb=msbm5 
\textfont\Bbbfam=\tenBbb 
\scriptfont\Bbbfam=\sevenBbb 
\scriptscriptfont\Bbbfam=\fiveBbb

\newcommand{\N}     {\mathbb{N}} 
\renewcommand{\P}   {\mathbb{P}}

\def\1{{\mathchoice {1\mskip-4mu\mathrm l}      % Blackboard bold 1 
{1\mskip-4mu\mathrm l} 
{1\mskip-4.5mu\mathrm l} {1\mskip-5mu\mathrm l}}} 
\newcommand{\ssup}[1] {{{\scriptscriptstyle{({#1}})}}} 
\def\comment#1{} 
\newtheoremstyle{thm}{2ex}{2ex}{\itshape\rmfamily}{} 
{\bfseries\rmfamily}{}{1.7ex}{} 
 
\newtheoremstyle{rem}{1.3ex}{1.3ex}{\rmfamily}{} 
{\itshape\rmfamily}{}{1.5ex}{}

\newtheorem{theorem}{Theorem}[section] 
\newtheorem{Lemma}[theorem]{Lemma} 
\newtheorem{proposition}[theorem] {Proposition} 
 
\newtheorem{remark}[theorem]  {Remark} 
\newtheorem{definition}[theorem] {Definition}

\theoremstyle{definition}

\renewcommand{\section}{\secdef\sct\sect} 
\newcommand{\sct}[2][default]{\refstepcounter{section} 
\vspace{0.8cm} 
\setcounter{equation}{0} 
\centerline{ %\large 
\large\scshape \arabic{section}.\ #1} 
\vspace{0.2cm}} 
\newcommand{\sect}[1]{ 
\vspace{0.8cm} 
\centerline{\large\scshape #1} 
\vspace{0.2cm}} 
 
\renewcommand{\subsection}{\secdef \subsct\sbsect} 
\newcommand{\subsct}[2][default]{\refstepcounter{subsection} 
\nopagebreak 
\vspace{0.5\baselineskip} 
{\flushleft\bf \arabic{section}.\arabic{subsection}~\bf #1  } 
\nopagebreak} 
\newcommand{\sbsect}[1]{\vspace{0.1cm}\noindent 
{\bf #1}\vspace{0.1cm}}

\renewcommand{\subsubsection}{% 
\secdef \subsubsect\sbsbsect} 
\newcommand{\subsubsect}[2][default]{% 
\refstepcounter{subsubsection} 
\nopagebreak 
\vspace{0.1\baselineskip} 
\nopagebreak 
{\flushleft 
\sffamily\slshape 
\arabic{section}.\arabic{subsection}.\arabic{subsubsection} 
\ % 
\sffamily #1\/.}\ } 
\newcommand{\sbsbsect}[1]{\vspace{0.1cm}\noindent 
{\bf #1}\ } 
 
\renewcommand{\d}{{\rm d}}

\newcommand{\Ccal}   {{\mathcal C }}

\newcommand{\Gcal}   {{\mathcal G }}

\newcommand{\df}      {{\,\stackrel{def}{=}\,}}
\begin{document}
\title[BOUNDS ON THE SPEED OF PROCESSES ON REGULAR TREES] {BOUNDS ON THE SPEED AND ON REGENERATION TIMES FOR CERTAIN PROCESSES ON REGULAR TREES}
\author[Andrea Collevecchio, Tom Schmitz]{}

\maketitle
\begin{center}
\thispagestyle{empty} 
\vspace{0.2cm} 
 
\centerline {\sc By  Andrea Collevecchio\footnote{Dipartimento di Matematica applicata, Universit\`a Ca' Foscari -- Venice, Italy.~~{\tt collevec@unive.it}}, Tom Schmitz\footnote{Max Planck Institute for Mathematics in the Sciences -- Leipzig, Germany.  {\tt schmitz@mis.mpg.de}}}

\vspace{0.4cm}

\centerline{\small(\version)}
\vspace{.5cm}
\end{center}
\[\]
\begin{abstract}
 We develop a  technique that  provides a lower bound on the speed of transient random walk in a random environment on regular trees. A refinement of this technique yields upper bounds on the first regeneration level and regeneration time. In particular, a lower and upper bound on the covariance in the annealed invariance principle follows. We emphasize the fact that our methods are general and also apply in the case of once-reinforced random walk.  Durrett, Kesten and Limic \cite{DKL2002} prove an upper bound of the form $b/(b+\delta)$ for the speed on the $b$-ary tree, where $\delta$ is the reinforcement parameter. For $\delta>1$ we provide a lower bound of the form $\gamma^2\,b/(b+\delta)$, where $\gamma$ is the survival probability of an associated branching process. 
 %In the latter case, the formula simplifies to $(1-\alpha)^{2}\frac b{b+\delta}$, where $\alpha$ is the smallest positive root of a known polynomial. 
 \end{abstract}

\vspace{.3cm}


\section{Introduction}
\label{intro}
\noindent  Random procesess with long memory have gained considerable attention in the recent past. Two emblematic examples of such processes are random walks in a random environment and reinforced processes. Although considerable progress has been achieved, there are many basic questions that remain open. %especially on $\mathbb Z^d$ in higher dimensions. 
We refer  to the overviews by Sznitman \cite{Szn2002} and Zeitouni \cite{Zeit2004},\cite{Zeit2006} for random walk in a random environment on $\mathbb Z^d$, and by Pemantle \cite{P2007} for reinforced processes on $\mathbb Z^d$ and on trees.\\
In this article we look at certain transient processes on regular trees, more precisely at random walk in a random environment, and at once-reinforced random walk. 
An important question is to obtain an explicit expression for the speed (if at all it exists), or at least to get good estimates. 
%In a first step one is of course interested in proving that the speed at all exists. This is often done with the help of regeneration techniques, and the speed is then given in terms of moments of regeneration levels and regeneration times, which are
This is in general a hard question, even for Markov chains as the biased random walk on a general tree, i.e.~a graph without cycles. For this model  there is in general no explicit expression for the speed, and often only an upper bound is at hand. It is in general hard to find a lower bound, and
% it is harder to derive a lower bound than an upper bound on the speed. 
we refer to Chen \cite{C2002} for several examples. We also point out to random walks on general graphs (Vir\'ag \cite{V2000}) where basically no lower bound on the speed is available.
%On $\mathbb Z^d$, this question is very hard for both these processes. Even worse, there is  no law of large numbers that holds in full generality (except for random walk in random environment when $d=2$).
%, and no recurrence-transience dichotomy. 

For random walk in a random environment, the speed is explicitly known only in one-dimensional models. On $\mathbb Z^d,\,d \ge 2$, not much is known about the speed, and even worse, if $d \ge 3$, it is still open if a law of large numbers with constant speed holds, see \cite{Szn2002},\cite{Zeit2004},\cite{Zeit2006}. On regular trees however, a law of large numbers holds, see \cite{G2004}, and transience implies that the speed is positive. This follows from  Theorem 1.1 in Aid\'ekon \cite{A2008} that treats the more general setting of Galton-Watson trees. 
One of our goals is to find  a {\it lower} bound on the speed for random walks in a random environment on regular trees. Our approach is general and we apply it to another class of processes with long memory: once edge-reinforced random walk.
%It is in this setting that we aim for a {\it lower} bound on the speed.
%We also refer the reader to Hu and Shi \cite{HS2007a},\cite{HS2007b} for an anlysis of the recurrent regime. 
Once edge-reinforced random walk on regular trees is transient, and has positive speed, see Theorem 1 and 2 in Durrett, Kesten and Limic \cite{DKL2002}. They propose an upper bound on the speed, but no lower bound that is always positive is at hand. With similar techniques than in the setting of random walk in random environment, we derive a lower bound. % that is expressed in terms of the survival probability of a related Galton-Watson process. 

%Our main task is to derive a lower bound on the speed for these two classes of processes. 
%Although random walk in random environment and once reinforced random walk are quite different, our methods apply in both settings. A similar idea already appeared in Collevecchio \cite{C2006b} for certain processes on Galton-Watson trees.
%A recurrence-transience dichotomy for biased random walk, transience for once-reinforced random walk and for vertex-reinforced jump process follow quite easily in \cite{C2006b}.
%An auxiliary  branching process is constructed to prove a recurrence-transience dichotomy for biased random walk, transience for once-reinforced random walk and for vertex-reinforced jump process. These results were already partially known, see \cite{C2006b} for further details, but the underlying idea is simple and either extends already known results or simplifies their proofs. 
%Our methods are inspired by ideas from Collevecchio \cite{C2006b}. 
%., which applies to  random walk in a random environment and to once-reinforced random walk. 
In order to provide a lower bound on the speed, it is instrumental to find a lower bound for the escape probability from the root, as well as an upper bound for the expected number of returns to the root. % see Proposition \ref{thm:1} and Theorem \ref{thm:ORRW}. 
Both these bounds are obtained with the help of an auxiliary branching process that already appeared in Collevecchio \cite{C2006b}. In particular the escape probability is bounded from below by the survival probability of the branching process, see the Propositions \ref{prop:branching} and \ref{lemma:weights}. For once-reinforced random walk, the branching process can be constructed in such a way that its survival probability is always positive, whereas for random walk in random environment we need additional assumptions.  

By a refinement of our methods, we are moreover able to derive a {\it common} explicit upper bound on all the moments of a first regeneration time $\tau_1$. These bounds are general and hold for  random walk in a random environment as well as for once edge-reinforced random walk, see Theorem \ref{thm:tau}. In words, this first regeneration time is the first time the height of the walk reaches a new maximum, and from then on never backtracks below this maximum. 
Regeneration times enjoy a wide-spread use in different settings, and we refer for instance to Lyons, Pemantle and Peres \cite{LPP1996} for biased random walk on a Galton-Watson tree, to Durrett, Kesten and Limic \cite{DKL2002} for once-reinforced random walk on a regular tree, and to Sznitman \cite{Szn2002} for random walk in a random environment on $\mathbb Z^d$.\\
The main step is to derive an explicit upper exponential tail on the first regeneration level $\ell_1$, defined as $\ell_1=|X_{\tau_1}|$, where $|\cdot|$ denotes the height of a vertex, see Theorem \ref{codaelle}. We inspire ourselves from Collevecchio \cite{C2009}, where a similar technique was introduced, although in the setting of the vertex-reinforced jump process. 
%The bound on the moment of $\tau_1
Let us mention that a detailed analysis of the tail behaviour of the first regeneration time is presented in Proposition 2.1 and 2.2 in Aid\'ekon \cite{A2009}, revealing an exponential and a subexponential regime on regular trees. We emphasize that we obtain {\it explicit} upper bounds on all moments of the first regeneration time under certain assumptions, in contrast to \cite{A2009}, where only the finiteness of the moments follows.
In particular, these bounds on the first regneration level resp.~regeneration time imply a lower and an upper bound on the covariance of the Brownian motion that appears as the limiting object in an annealed invariance principle, see Theorem \ref{thm:3} and Proposition \ref{prop:5}.
 %On these graphs, the genuinly irreversible random walk in a random environment on $\mathbb Z^d,\,d \ge 2$, becomes reversible.
%In the setting of random walk in a random environment on $\mathbb Z^d$, there is for instance no recurrence-transience dichotomy in dimension two or higher, and there is no strong law of large numbers in dimension three or higher. In particular, if the so-called condition $(T)$ holds in dimension two or higher, then it is known that a law of large numbers with non-vanishing speed holds, see \cite{Szn2002}. However not much more can be said about the speed, and even worse, the speed can have opposite direction to the expected local drift, see Bolthausen, Sznitman and Zeitouni \cite{BSZ2003}. In contrast, these questions have been solved in the one-dimensional case, and the speed is explicitly known.\\
 
This article is organised as follows. In Section 2, we provide a lower bound on the speed for random walk in a random environment and for once edge-reinforced random walk, and in Section 3 we derive moment bounds on the first regeneration time that are completely general and hold for random walk in a random environment and for once edge-reinforced random walk.
%\newpage
%%%%%%%%%%%%%%%%%%%%%%%%%%%%%%%%%%%%%%%%%%%%%%%%%%%%%%%%%%%%%%%%%%%%%%%%%%%%%%%%%%%%%%%%%%%%%%%%%%%%%%%%%%%%%%%%%%%%%%%%%%%%%%%%%%%%%%%%%%%%%%%%%%%%%%%%%%%%%%%%%%%%%%%%%%%%%%%%%%%%%%%%%%%%%%%%%%%%%%%%%%%%%%%%%%%%%%%%%%%%%%%%%%%%%%%%%%%%%%%%%%%%%%%%5

\section{On the speed}

Let us start by introducing some notation. 
Consider the $b$-ary regular tree $\Gcal_b$ with root $\rho$. %We suppose that the root has an additional parent $\overleftarrow{\rho}$.
We assume that the root $\rho$ has a parent $\overleftarrow{\rho}$. Hence each vertex in the tree is connected to $b+1$ vertices, except for $\overleftarrow{\rho}$, that is only connected to $\rho$.
%Each vertex has exactly $b+1$ neighbors with the exception of $\overleftarrow{\rho}$  which is connected  only to  $\rho$. 
%We define the distance between two vertices in $\Gcal_{b}$  as the number of edges on the unique self-avoiding path connecting them. 
For any vertex $\nu$, denote by $|\nu|$ its distance to the root, i.e.~the number of edges on the unique self-avoiding path connecting $\nu$ and $\rho$. Level $i$ is  the set of vertices $\nu$ such that $|\nu|=i$, with the exception that $|\overleftarrow \rho|=-1$.
 For $\nu \neq \overleftarrow \rho$, define $\overleftarrow{\nu}$, called the parent of $\nu$,  to be the unique vertex  at level  $|\nu|-1$ connected to $\nu$. We say that $\nu$ is a child of $\overleftarrow{\nu}$. 
 % We denote by $\stackrel{\Leftarrow}{\nu}$ the grandfather of $\nu$, that is the parent of its parent. 
 %Moreover, we choose a complete order $\preceq$ among the vertices of $\Gcal_{b}$. For any vertex $\nu \neq \rho$ we use this order to label its children by $\overrightarrow{\nu}^{\ssup 1},\overrightarrow{\nu}^{\ssup 2}, \ldots, \overrightarrow{\nu}^{\ssup b}$ , with 
 %$ \overrightarrow{\nu}^{\ssup i} \preceq \overrightarrow{\nu}^{\ssup j}$ if  $i \le j$. 
We say that a vertex $ \nu_{0}$ is a  descendant of the vertex $\nu$ if the latter lies on the unique self-avoiding path connecting $\nu_{0}$ to $\rho$, and $ \nu_{0} \neq \nu$. In this case, $\nu$ is said to be an ancestor of $\nu_{0}$.
For any vertex $\mu$, let $\L_{\mu}$  be the subtree of  $\Gcal_{b}$  consisting of $\mu$,  its descendants and the edges connecting them, i.e.~the $b$-ary subtree rooted at $\mu$. Let  $\overleftarrow{\L}_{\mu}$ be the smallest subtree of $\Gcal_{b}$ containing $\L_{\mu}$ and the vertex  $\overleftarrow{\mu}$. 

\subsection{Random Walk in Random Environment}
Let us define the random environment. To each vertex $\nu$, different from $\overleftarrow \rho$, we assign a $b$-dimensional random vector with positive entries
 $$
  \mathbf{A}_{ \nu} \df ( A^{\ssup 1}_{\nu},A^{\ssup 2}_{\nu},\ldots A^{\ssup b}_{\nu}).
 $$ 
 We assume that these vectors  are i.i.d.~~under the measure $\P$.
Moreover, following Lyons and Pemantle \cite{LP1992}, we assume that the coordinates are identically distributed.  The random environment $\omega$ is defined by $\omega(\overleftarrow \rho, \rho)=1$ and for any vertex $\nu \neq \overleftarrow \rho$,
\begin{equation}
\label{eq:1.1}
\omega(\nu,\overrightarrow \nu^{\ssup i})=\frac{A^{\ssup i}_{\nu}}{1+\sum_{j} A^{\ssup j}_{\nu}}; \quad \omega(\nu,\overleftarrow \nu)=\frac 1{1+\sum_{j} A^{\ssup j}_{\nu}}.
\end{equation}
%It is often convenient to use the short-hand notations
%\begin{equation}
%\label{eq:1.2}
%$\omega_\nu \df \omega(\overleftarrow \nu,\nu)$, and for a generic copy of $\omega_\nu$,  we write $\omega_{\text{up}}$.
%\end{equation}
For a vertex $\nu$ we define the Markov chain $\{ X_{n}, n \ge 0\}$ started at $\nu$ by
\begin{eqnarray*}
\mathbf{P}_{\nu, \omega}(X_{0} = \nu) &=& 1\\
\mathbf{P}_{\nu, \omega}(X_{n+1} = \mu_{1} | X_{n} = \mu_{0}) &=& \omega(\mu_{0}, \mu_{1}), 
\end{eqnarray*}
  for any pair of neighbors $\mu_{0}, \mu_{1}$. We introduce further the annealed measure as the semi-direct product $\mathbf{P}_{\nu}=\mathbb P \times \mathbf{P}_{\nu,\omega}$.  We write $\mathbf{P}_{\omega}$ and $\mathbf{P}$ for $\mathbf{P}_{\rho,\omega}$ resp.~$\mathbf{P}_{\rho}$. We also write  
$A$ and $\mathbf A=(A^{(1)},\ldots,A^{(b)})$ for a generic copy of  $A^{\ssup i}_{\nu},\,1 \le i \le b$, respectively for a generic copy of $\mathbf A_\nu=(A^{\ssup 1}_{\nu},\ldots,A^{\ssup b}_{\nu})$.
%\begin{theorem}\label{stronglaw} 
%Let $\mathbf{X}$ be RWRE on the binary tree.  Suppose that 
%\begin{equation}\label{momcond}
%\beta \df \mathbb{E}\Big[\frac{\xi^{\ssup \nu}_{0}}{ \xi^{\ssup \nu}_{1}}\Big] <2.
%\end{equation}
%Then there exists a  constant $  K \in (0, \infty)$ such that   
%\begin{eqnarray}
%\label{slln} \frac{|X_k|}{k} \rightarrow K  \qquad \mbox{a.s.}.
%\end{eqnarray}
%\end{theorem}
We introduce the  hitting  times of a vertex $\nu$ 
respectively of a level $i$
\begin{equation}\label{T}
 T(\nu) \df \inf \{ k \ge 0 \colon X_{k} = \nu\} \qquad \mbox{and} \qquad  T_i \df \inf \{ k \ge 0 \colon |X_k| = i \}.
\end{equation}
We further introduce the respective return times
\begin{equation}
\label{beta0}
 D \df \inf\{ n \ge 1 \colon X_{n} = \overleftarrow{X_0}\}, \qquad D(\nu) \df \inf\{ n \ge 1 \colon X_{n-1}=\nu,\,X_n=\overleftarrow \nu \},
\end{equation}
and the annealed return probability
\begin{equation}
\label{beta}
\quad \beta \df \mathbf{P}(D < \infty).
\end{equation}
To each ordered pair of neighbors $ \nu,\mu \in $Vert$(\Gcal_{b})$ assign 
%a Poisson process $P(\nu,\mu)$ of rate 1. 
a collection of independent exponentials $h_{k}(\nu,\mu)$, $ k \ge 0$, each with mean one.
We assume that all   these collections are independent. 
%Call $h_{k}(\nu,\mu)$, with $ k \ge 0$, the inter-arrival times of $ P(\nu,\mu)$.
Using these exponentials, we now provide a construction of random walk in random environment on an arbitrary subtree (see \cite{S1994}  for a similar construction for reinforced processes). 
\begin{definition}\rm
\label{ext}
(Extension $\mathbf{Y}^{\Ccal}$)
Fix a subtree $\Ccal$ of $\Gcal_b$. The extension $\mathbf{Y}^{\Ccal}$ of $\mathbf{X}$ on the subtree $\Ccal$ is defined as follows. 
Fix a starting point $\eta$ in $\Ccal$, i.e.~$Y^\Ccal_0=\eta$. We define $\mathbf Y^\Ccal$ iteratively in the following way.
 Let $ s_{1}(\nu) $ be the first time $\mathbf{Y}^{\Ccal}$
  reaches some vertex $\nu$. Define $ N_\nu^\Ccal$ %=\{\overleftarrow \nu, \overrightarrow{\nu}^{\ssup i},i=1,\ldots,b\}$ 
to be the set of neighbors of $\nu$ in $\Ccal$.
 The first jump after $s_{1}(\nu)$  is towards the neighbor $\mu \in N_{\nu}^\Ccal$ for which the following minimum
%Let $\Ccal$ be a subtree rooted at $\mu$. The extension $ \mathbf{X}^{\Ccal}$ starts at $\mu$ and is defined recursively as follows. Fix  a vertex $\nu$ of $\Ccal$. Let $ s_{1} $ be the first time $\mathbf{X}^{\Ccal}$
%  reaches $\nu$. Define $ N_\nu=\{\overleftarrow \nu, \overrightarrow{\nu}^{\ssup i},i=1,\ldots,b\}$ to be the set of neighbors of $\nu$ in $\Ccal_{b}$.
% The first jump after $s_{1}$  is towards the neighbor $\eta \in N_\nu \cap $Vert$(\Ccal)$  for which the following minimum
\begin{equation}
\label{eq:expjump1}
\min_{\eta \in N_\nu^\Ccal} \frac{h_{1}(\nu,\eta)}{\omega(\nu,\eta)}
%h_{1}(\nu,\overleftarrow{\nu}) \wedge \min_{i}  \big(h_{1}(\nu,\overrightarrow{\nu}^{\ssup i})/A^{\ssup i}_{\nu}\big),
\end{equation}
%where $a \wedge b= \min\{a,b\}$. 
is a.s.~attained.
We define  $s_{k}(\nu),\,\, k \ge 2$, inductively via
%Suppose we defined $s_{j}, 1 \le j \le k-1$, and let 
\begin{eqnarray*}
 && s_{k} \df  \inf \big\{ n > s_{k-1} \colon Y^{\Ccal}_{n} = \nu \big\}, \mbox{ and} \\
 &&  j_{k}(\nu,\mu)  \df 1+ \mbox{ number of times $\mathbf{Y}^{\Ccal}$ jumped from $\nu$ to its  neighbor $\mu$  by time $s_{k}$}.
 \end{eqnarray*}
The first jump after $s_{k}$  is towards the neighbor $\mu$ for which the following minimum 
\begin{equation}
\label{eq:expjump2}
\min_{\mu \in N_\nu^\Ccal} \frac{h_{j_{k}}(\nu,\mu)}{\omega(\nu,\mu)}
%h_{j(\nu,\overleftarrow{\nu})}(\nu,\overleftarrow{\nu})\wedge \min_{i} \big( h_{j(\nu, \overrightarrow{\nu}^{\ssup i})}(\nu,\overrightarrow{\nu}^{\ssup i})/A^{\ssup i}_{\nu} \big)
\end{equation}
is a.s.~attained. With a slight abuse of notation, we denote the quenched and annealed law of the extension $\mathbf Y^\Ccal$ again by $\mathbf P_{\cdot,\omega}$ resp.~$\mathbf P_\cdot$.
\qed \end{definition} 
%Denote with $N_{\mu}$  the set of neighbors of the vertex $\mu$ in $\Gcal_b$. 
%Denote the children of a vertex $\nu$ with $\nu^{(i)}$. 
\begin{remark}\rm
\label{extcoupl} 
The extension processes will play a crucial role in our proofs. They are coupled to the original process $\mathbf{X}$ 
in the following sense.
%because their jumps are generated by the same collection of exponentials. To be more precise, 
Let  $\mathbf{Y}^{\Ccal}$ be the extension of $\mathbf{X}$ on $\Ccal$, started at a vertex $\nu$ in $\Ccal$. Denote with $\theta_\cdot$ the canonical time shift, and suppose that $\mathbf{X}$ hits $\nu$. Since both processes are generated by the same exponential variables, it follows that $\mathbf{Y}^{\Ccal}$ coincides with the process $\mathbf{X} \circ \theta_{T(\nu)}$, of course only observed on the subtree $\Ccal$, which is called {\it restriction } process. For a rigorous definition of restriction process see \cite{C2006} or \cite{DV2002}. Extension processes were used  in \cite{C2009} to prove the strong law of large numbers for vertex jump-reinforced processes. 
%This coupling lasts up to the last visit of  $ \mathbf{X}$ to this   subtree.  
\qed 
\end{remark} 
A child $\nu^{(j)}$ of $\nu$ is called a {\it first child}  if 
it is a.s.~the minimiser of 
\begin{equation}
\label{firstch}
\min_{1 \le i  \le b}\frac{h_{1}(\nu,\nu^{(i)})}{\omega(\nu,\nu^{(i)})}\quad a.s.
\end{equation}

Let us now turn to the lower bound on the speed.
Lyons and Pemantle \cite{LP1992} (see also Menshikov and Petritis \cite{MP2002}) established the following recurrence-transience dichotomy:
\begin{equation}
\mathbf{X} \text{ is transient if $\inf_{0 \le t \le 1}\mathbb E[A^t]>\tfrac 1b$, and recurrent otherwise.}
\end{equation}
Our standing assumption is that the walk is transient. Gross \cite{G2004} proves a strong law of large numbers 
\begin{equation}
v \df \lim_{n \to \infty}\frac{|X_n|}n \ge 0 \quad \mathbf P-\text{a.s.}
\end{equation}
The natural question to ask now is in which cases $v$ is positive. This question was answered recently in Aid\'ekon \cite{A2008} in the more general setting of Galton-Watson trees. In our setting, on regular trees, it turns out that $v$ is always positive, see Theorem 1.1 in \cite{A2008}. 
%Under certain conditions to be formulated below, 
We will now derive a lower bound on the speed $v$.
For $n\ge 1$, we define
\begin{equation}
\label{eq:L}
L(\nu, n) \df \sum_{j=0}^{n} \1_{\{X_{j} = \nu\}}, \quad \text{and} \quad L(\nu) \df \sum_{j=0}^{\infty} \1_{\{X_{j} = \nu\}},
\end{equation}
the number of visits to $\nu$ by time $n$, resp.~the total number of visits.
Under transience, it is well-known that $v= \lim_{n \ti} |X_{n}|/n$ exists. Here is the main result of this subsection.
\begin{proposition} 
\label{thm:1}
Under transience, it holds that
%that $\mathbb{E}[A^{-2}] < \infty$. 
\begin{equation} 
v   \ge \frac{1-\beta}{\mathbf E[L(\rho)]}>0\,\quad \mathbf P-\text{a.s.}
%\frac{1- \beta}{\theta(2,1)^{1/2}(1-\beta)^{b/2}+
%      \sum_{n=2}^\infty \theta(2,n)^{1/2}\,\beta^{b^{n-2}/2}(\sum_{i=1}^b(-1)^{i-1} {b \choose i}\beta^{b^{n-2}(i-1)} 
%      )^{1/2}}>0.
\end{equation}
%\begin{equation}\label{Theta}
%\begin{aligned}
%\Theta(p)=&\sum_{\heap{F \in \Pcal_{b}}{F \neq \emptyset}}  c(p) \mathbb{E}\Big[\big(\sum_{\nu \in F}\omega(\rho,\nu) \big)^{-p}\Big] \beta^{b - |F|} (1- \beta)^{|F|}\\ &+  \big(c(p) \mathbb{E}[\omega^{-p}(\rho, \overrightarrow{\rho}^{\ssup 1})] \big)^{1/2} \Big(\sum_{k=2}^{\infty}  \Big(1 + c(p)k^{p-1}\frac{\mathbb{E}[A^{-p}]^{k} -1}{\mathbb{E}[A^{-p}] -1}\Big)^{1/2}   \beta^{b^{k-1}/2}\Big),
%\end{aligned}
%\end{equation}
%where $\Pcal_{b}$ denotes the power set of $\{\overrightarrow{\rho}^{\ssup 1},\ldots, \overrightarrow{\rho}^{\ssup b}\}$, and $|F|$ is the cardinality of the set $F$.
%\sum_{k=1}^{\infty} \big(\Theta(p,k)\big)^{1/p} \beta^{(p-1)b^{k-2}/p}
\end{proposition}
%Recall $\beta$ in (\ref{beta}). 
%The  proof of  Theorem \ref{thm:1} relies on the following  key result.
%We postpone the proof of Proposition \ref{prop:momtime}, and provide first the proof of Theorem \ref{thm:1}.
%In the proof of Theorem \ref{thm:1}, we apply this proposition with $p=1$. Later on, we will use it for general $p \ge 1$.
\noindent Before proving Proposition \ref{thm:1}, we provide first a lemma. Let 
\begin{equation}
\label{eq:Pi}
\Pi_{k}=\sum_{\nu \colon |\nu|=k}\1_{\{T(\nu)<\infty\}}
\end{equation}
be the number of vertices visited at level $k$. Recall $\beta$ in (\ref{beta}). We have 
\begin{Lemma}
\label{lemma:A.1}
Assume transience, i.e.~$\beta<1$. Then $\Pi_{k}$ is stochastically dominated by a geometric random variable with parameter $1-\beta$.
\end{Lemma}
\begin{proofsect}{Proof}
One vertex at level $k$ is visited for sure. Call this vertex $\sigma_1$.  Notice that, after $T(\sigma_1)$, a necessary condition to visit a further vertex at level $k$ is that the walk returns to the parent of $\sigma_1$.
To obtain an upper bound for $\Pi_{k}$, we adopt the following strategy. If the walk returns to the parent of $\sigma_1$, we consider the extension $Y^{(\sigma_1)}$ of $\mathbf X$ to the subtree obtained by cutting the subtree $\Lambda_{\sigma_1}$. 
This ensures that the second visit at level $k$ will be at a new vertex $\sigma_2$, different from $\sigma_1$.
We repeat this procedure iteratively, and it clearly yields an upper bound on the number of vertices $\sigma_i$ visited at level $k$.
Each time a new vertex $\sigma_i$ is visited, there is a chance of escape to infinity with annealed probability $1-\beta>0$, because of stationarity.
Since all subtrees $\Lambda_{\sigma_i}$ are disjoint, the trials of escape are independent. It follows  that $\Pi_k$ is dominated by a geometric with parameter $1-\beta$. This ends the proof.
\qed \end{proofsect}

\begin{proofsect}{Proof of Proposition~~\ref{thm:1}} 
Notice that 
\begin{equation}
\label{eq:1}
\lim_{n \ti} \frac{T_{n}}n = 1/v \qquad \mathbf{P} \mbox{ a.s.}
\end{equation}
Label the vertices at level $k$  by $\nu_{k,1}, \nu_{k,2}, \ldots,  \nu_{k,b^{k}}$. We have that for $n \ge 1$,
\begin{equation}
\label{eq:0.2}
\begin{aligned}
\mathbf{E}[T_{n}] \le  1+\mathbf{E}[L(\overleftarrow{\rho})]+ \sum_{k=0}^{n-1} \sum_{j=1}^{b^{k}} \mathbf{E}\big[ L(\nu_{k,j}) \1_{\{T(\nu_{k,j})<\infty\}} \big].
\end{aligned}
\end{equation}
Fix a vertex $\nu$, and define $\widetilde L(\nu)$ to be the total time spent in the vertex $\nu$ by the extension of $\mathbf X$ to $\overleftarrow{\Lambda}_{\nu}$ started at $\nu$. Then $L(\nu) \le \widetilde L(\nu)$, and the law of 
$\widetilde L(\nu)$ under $\mathbf P_\nu$ is equal to the law of $L(\nu)$ under $\mathbf P$.
Moreover the random variables $\widetilde L(\nu)$ and $\1_{\{T(\nu)<\infty\}}$ are independent under the annealed measure. We use independence, and then stationarity, and obtain that the sum on the right-hand side of \eqref{eq:2} is smaller than
\begin{equation}
\label{eq:0.3}
\begin{aligned}
 \sum_{k=0}^{n-1} \sum_{j=1}^{b^{k}} \mathbf{E}[\widetilde{L}(\nu_{k,j})]\,\mathbf{P}\big[T(\nu_{k,j})<\infty \big]
= \mathbf{E}[L(\rho)]\,\sum_{k=0}^{n-1}\mathbf{E}[\Pi_{k}]
\le  \mathbf{E}[L(\rho)]\frac {n}{1-\beta},
\end{aligned}
\end{equation}
where in the last step we used Lemma \ref{lemma:A.1}. Using \eqref{eq:0.2} and \eqref{eq:0.3}, and by Fatou's lemma, we obtain that $\mathbf P$-a.s.,
\begin{equation}
 \lim_{n \to \infty} T_n/n \le \liminf_{n \to \infty}\mathbf{E}[T_n/n]=\mathbf{E}[L(\rho)](1-\beta)^{-1}.
\end{equation}
The claim of the theorem follows now from (\ref{eq:1}).
 \qed
\end{proofsect}
Our main task is now to derive upper bounds on $\beta$ and on the expectation of $L(\rho)$.

%%%%%%%%%%%%%%%%%%%%%%%%%%%%%%%%%%%%%%%%%%%%%%%%%%%%%%%%%%%%%%%%%%%%%%%%%%%%%%%%%%%%%%%%%%%%%%%%%%%%%%%%%%%%%%%%%%%%%%%%%%%%%%
\subsubsection{Estimates on the return probability $\beta$}
\noindent In the last section, we provided a lower bound in terms of the annealed return probability $\beta$. In this section, we will derive an upper bound on $\beta$ in terms of the extinction probability $\alpha$ of a certain branching process, in the spirit of Collevecchio \cite{C2006b}. This allows to obtain an explicit  lower bound on the speed. \\
Let us start by constructing the branching process. 
\begin{definition}(Color scheme)\rm
\label{def:color}
Fix an integer $\psi \ge 1$, and denote with $\mathbf Y(\nu,\mu)$ the extension of $\mathbf X$ to the unique ray connecting the vertices $\mu$ and $\nu$.
%For any vertex $\nu$, with $|\nu| \ge 1$,  define $H_{\nu}:= \inf\{ n \ge t(\nu) \colon X_{n} = {\rm par}(\nu)\}$. 
We introduce the following color scheme. A vertex $\nu$ at level $\psi$ is colored if and only if the $\mathbf Y(\overleftarrow \rho,\nu)$, started at $\rho$, hits $\nu$ before $\overleftarrow \rho$. A vertex $\nu$ at level $k\psi,\,k \ge 2$, is colored if and only if
\begin{itemize}
\item its ancestor at level $(k-1)\psi$, say $\mu$, is colored, and 
\item  $Y(\overleftarrow \mu,\nu)$, started at $\mu$, hits $\nu$ before  $\overleftarrow \mu$.
\end{itemize}
All the other vertices are uncolored, and only vertices that are at a level $k\psi$, $k\ge 1$, can be colored. 
\qed \end{definition}
Under the {\it annealed} measure, the number of colored vertices form a homogeneous branching process, since the offspring is each time determined by disjoint parts of the environment. We denote this branching process with $\mathbf Z_\psi$. 
%Later on, we will specify assumptions under which the branching process $\mathbf Z$ is supercritical. 
%Assume for the moment that $\mathbf Z$ survives. 
%Notice that all vertices $X_{T_{k\psi}},\,k \ge 1,$ are colored. In particular, if the branching process survives, then each level $k\psi,\,k \ge 1$, is hit before returning to the parent of the root. In other words, $D=\infty$ $\mathbf P$-a.s., and the process $\mathbf X$ is transient. If we denote with $\alpha$ the extinction probability of $\mathbf Z$, then it follows that
We formulate the following
\begin{proposition}
\label{prop:branching}
Denote with $\alpha_\psi$ the extinction probability of $Z_\psi$. Then
%\begin{equation}
%\label{eq:alpha}
 $\beta \le \alpha_\psi$.
%\end{equation}
If moreover $\mathbb E[A^{-1}]<b$, then there is an integer $\psi \ge 1$ such that $\alpha_\psi<1$.
\end{proposition}
\begin{proofsect}{Proof}
Let us show that $\beta \le \alpha_\psi$ in the case $\alpha_\psi<1$ (otherwise there is nothing to prove). Assume that $\mathbf Z_\psi$ survives. 
Choose vertices $\mu$ and $\nu$ as in definition \ref{def:color}. 
By remark \ref{extcoupl}, the processes  $Y(\overleftarrow \mu,\nu)$  and $\mathbf{X}$ coincide, from the time the latter hits $\mu$ until its last visit to the path connecting $\overleftarrow \mu$ to $\nu$. It follows that, if $\mathbf X$ hits $\nu$ before $\mu$, then so does $Y(\overleftarrow \mu,\nu)$. 
It follows that all vertices $X_{T_{k\psi}},\,k \ge 1,$ are colored. 
In particular, if the branching process survives, then each level $k\psi,\,k \ge 1$, is hit before returning to the parent of the root. %In other words, $D=\infty$ $\mathbf P$-a.s. 
Hence $\{Z_\psi \text{ survives}\} \subseteq \{D=\infty\}$, and $\beta \le \alpha_\psi$ follows.
%Let us now compute the extinction probability $\alpha_\psi$. 
Let us now show that if $\mathbb E[A^{-1}]<b$, then we can find $\psi \ge 1$ such that $\mathbf Z_\psi$ is supercritical.
We choose a vertex $\mu$, and then a vertex $\nu$ at level $|\mu|+\psi$. Then the extension $\mathbf Y(\overleftarrow \mu,\nu)$, started at $\mu$, hits $\nu$ before $\overleftarrow \mu$ with (annealed) probability
\begin{equation}
\label{eq:ruin}
\mathbb{E}\big[ \big( \sum_{r=1}^{\psi+1} \prod_{j=1}^{r-1} A_{j}^{-1} \big)^{-1}\big],
\end{equation}
where $A_j, 1 \le j \le \psi$, is an enumeration of the variables $A$ along the ray connecting $\mu$ to $\nu$.
By Jensen's inequality,
\begin{equation}
\label{eq:2}
 \mathbb{E}\big[ \big( \sum_{r=1}^{\psi+1} \prod_{j=1}^{r-1} A_{j}^{-1} \big)^{-1}\big]
\ge \mathbb{E}\big[ \sum_{r=1}^{\psi+1} \prod_{j=1}^{r-1} A_{j}^{-1}\big]^{-1}.
\end{equation}
By independence,  we find for large $\psi$,
\begin{equation}
\label{eq:3}
 \mathbb{E}\big[ \sum_{r=1}^{\psi+1} \prod_{j=1}^{r-1} A_{j}^{-1}\big]
=\sum_{r=1}^{\psi+1}\mathbb E[A^{-1}]^{r-1} =\frac{1-\mathbb E[A^{-1}]^{\psi+1}}{1-\mathbb E[A^{-1}]}.
\end{equation}
By the assumption $\mathbb E[A^{-1}]<b$, we find that
\begin{equation}
\label{eq:4}
 \lim_{\psi \to \infty} b^{-\psi}\,\mathbb{E}\big[ \sum_{r=1}^{\psi+1} \prod_{j=1}^{r-1} A_{j}^{-1}\big]=0\,.
\end{equation}
Hence, if we choose $\psi$ large enough, then we can make sure that
\begin{equation}
\label{eq:5}
 b^{\psi}\,\mathbb{E}\big[ \big( \sum_{r=1}^{\psi+1} \prod_{j=1}^{r-1} A_{j}^{-1} \big)^{-1}\big]>1\,.
\end{equation}
Notice that the left-hand side of the last display is the expected offspring of the branching process $Z_\psi$, so that we can choose $\psi$ s.t.~$\mathbf Z_\psi$ is supercritical.  This finishes the proof of the proposition.
\qed \end{proofsect}
\begin{definition}\rm
%From now on, we use $\psi$  to denote the smallest positive integer satisfying (\ref{eq:5}), and 
We denote with
$\mathbf{p}:=\{p_{k},\, k \in \{0,1, \ldots, b^{\psi}\} \}$ the offspring distribution of the branching process $\mathbf Z_\psi$. The mean offspring is 
 \begin{equation}
\label{eq:m}
 m_\psi \df \sum_{k=0}^{b^{\psi}} k p_{k}.
 \end{equation}
\qed \end{definition}
Proposition \ref{prop:branching} implies that if $\mathbb E[A^{-1}]<b$, then there is $\psi \ge 1$ such that $m_\psi>1$.\\

%We use $\psi$  to denote the smallest positive integer satisfying \eqref{niceint}. Moreover we denote by 
%$\mathbf{p}:=\{p_{k},\, k \in \{0,1, \ldots, b^{\psi}\} \}$ the offspring distribution of the branching process of the colored vertices. We assume that $\mathbb{E}[\frac 1 A]< b$ holds, hence 
% \begin{equation}
% \label{eq:m}
% m \df \sum_{k=0}^{b^{\psi}} k p_{k}>1.
% \end{equation}
%\begin{proposition}\label{alpha} Set $\alpha$ to be the smallest positive solution of the equation
%\begin{equation}
%x= \sum_{i =0}^{b^{\psi}} x^{i} p_{i}.
%\end{equation}
%Then $\beta \le \alpha$.
%\end{proposition}
%%%%%%%%%%%%%%%%%%%%%%%%%%%%%%%%%%%%%%%%%%%%%%%%%%%%%%%%%%%%%%%%%%%%%%%%%%%%%%%%%%%%%%%%%%%%%%%%%%%%%%%%%%%%%%%%%%%%%%%%%%%
\subsubsection{An explicit upper bound on the expectation of $L(\rho)$}
Our standing assumption in the remaining subsections is that
\begin{equation}
\label{eq:A}
 \text{we can find $\psi \ge 1$ such that $\alpha_\psi<1$},
\end{equation}
where we recall $\alpha_\psi$ in Proposition \ref{prop:branching}. Condition (\ref{eq:A}) is in particular satisfied if 
$\mathbb E[A^{-1}]<b$. For $p \ge 1, n \ge 1$, we introduce the function
\begin{equation}
 \label{eq:theta}
 \theta(p,n) \df  
\begin{cases}%{ll}
c_p\,b\,\mathbb{E}[(1+\tfrac 1{\sum_{i=1}^b A^{(i)}})^p]\, n^{p-1}\,\frac{\mathbb{E}[A^{-p}]^{n} -1}{\mathbb{E}[A^{-p}] -1},&\text{ if }n \ge 2,\\
c_p\,\mathbb{E}[(1+\tfrac 1{\sum_{i=1}^b A^{(i)}})^p],&\text{ if }n=1,\\
\end{cases}
\end{equation}
where the r.h.s.~is infinite if $\mathbb{E}[A^{-p}]=\infty$, and the constants $c_p$ are introduced in Lemma \ref{lemma:A.2} in the Appendix. We have the following
\begin{proposition}
\label{prop:momtime} 
If $\mathbb {E}[A^{-p-\varepsilon}] < \infty$ for some $p \ge 1$ and some $\varepsilon >0$, then for all $n \ge 1$, $\theta(p+\varepsilon,n)<\infty$, and
\begin{equation}
\mathbf{E}[L(\rho)^{p}] \le \theta(p+\varepsilon,1)^{1/q}
%(1-\beta)^{b/q'}
+\,\sum_{n=2}^\infty \,\theta(p+\varepsilon,n)^{1/q}\alpha_\psi^{b^{n-2}/q'}\,
\Big(\sum_{i=1}^b (-1)^{i-1} {b \choose i} \alpha_\psi^{b^{n-2}(i-1)}\Big)^{1/q'},
\end{equation}
where $q=1+\varepsilon/p$, and $q'=1+p/\varepsilon$ is the dual of $q$.
\end{proposition}
%\begin{remark} 
%\label{rem:momtime}
%\rm
%In the proof of Proposition \ref{prop:momtime} below we use Cauchy-Schwarz's inequality in \eqref{eq:3.7}. If  instead we apply H\"older's inequality we can relax the assumption $\mathbb E[A^{-2p}]<\infty$ of Proposition \ref{prop:momtime} and replace it by the assumption $\mathbb E[A^{-p-\varepsilon}]<\infty$ for some $\varepsilon >0$.
%\end{remark}
%\subsection{Proof of Proposition~~\ref{momtime}.}
%Notice that condition \eqref{firstc} does not imply that the vertex fc($\nu$) is visited by the process. If $\mathbf{X}$ visits it, then it is the first among the children of $\nu$ to be visited.\\
Before proving Proposition \ref{prop:momtime}, we formulate an auxiliary result. We first introduce some notation. Fix $n \ge 2$. 
Choose $b$ distinct vertices $\nu_i,\,1 \le i \le b,$ at level $n$, with different ancestors at level one. More precisely, we choose $\nu_i$ with ancestor $\overrightarrow{\rho_i}$ at level one, and call this set of vertices $\mathcal A_n$. 
%From each child $\overrightarrow{\rho_{i}}$ of the root, fix a ray from $\overrightarrow{\rho_{i}}$ to level $n$, and 
We label the vertices on the ray connecting $\overrightarrow{\rho_i}$ to $\nu_i$ by $\sigma^{(i)}_j,\, 1 \le j \le n$, with $\sigma^{(i)}_1=\overrightarrow \rho_{i}$ and $\sigma^{(i)}_n=\nu_i$. Denote with $\Gamma_n$ the subtree composed by the root $\rho$, its parent $\overleftarrow \rho$, the vertices $\sigma^{(i)}_j, 1 \le j \le n, 1 \le i \le b$, and the edges connecting them. 
%$\mathcal A_n=\{\sigma^{(i)}_n, 1 \le i \le b\}$ denotes the set of vertices in $\Gamma_n$ at level $n$. 
For $n=1$, $\Gamma_1$ is simply the subtree composed by the root and its children, with the edges connecting them, 
and $\mathcal A_1$ is the set of children of the root. 
We denote with $\mathbf Y$ the extension of $\mathbf X$ to $\Gamma_n$, and we introduce
$\widetilde T_{\mathcal A_n}=\inf\{n \ge 0 \colon Y_n \in \mathcal A_n\}$, and $\widetilde{T}(\rho) \df \inf\{ n\ge 1 \colon Y_{n} =\rho\}$. We further define
$$
\widetilde{L}(\rho, \widetilde T_{\mathcal A_n})  \df \sum_{i=0}^{\infty} \1_{\{Y_{i} = \rho,\,i<\widetilde T_{\mathcal A_n}\}}\,.
$$
Recall $\theta(p,n)$ in \eqref{eq:theta}. We have the following
\begin{proposition}
\label{prop:3.1}
If $\mathbb{E}[ A^{-p}]< \infty$ 
%and that $\mathbb{E}[\omega^{-p}_{\text{up}}]<\infty$ 
for some $p\ge 1$, then 
\begin{equation}
\label{theta}
\mathbf{E}\big[\widetilde{L}(\rho, \widetilde T_{\mathcal A_n})^{p}\big] \le \theta(p,n)<\infty\,.
%b\,c_p\, \theta(p,n)=\,c_p\, b\, \mathbb{E}[(1+\tfrac 1{\sum_{i=1}^b A^{(i)}})^p]\, n^{p-1}\,\frac{\mathbb{E}[A^{-p}]^{n} -1}{\mathbb{E}[A^{-p}] -1},
\end{equation}
%and for $n =1$, we have 
%\begin{equation}
%\label{theta1}
%\mathbf{E}\big[\widetilde L \big(\rho,\widetilde T_{\mathcal A_n} \big)^{p}\big] \le c_p\,\theta(p,n))=c_p \,\mathbb{E}[(1+\tfrac 1{\sum_{i=1}^b A^{(i)}})^p].
%\end{equation}
\end{proposition}
\begin{proofsect}{Proof of Proposition \ref{prop:3.1}} Fix $n \ge 2$. 
%Each time $\mathbf{Y}$ hits the root, it has 
To escape from the root, the walk $\mathbf Y$ has to jump to one of the children of the root, and then hit the set $\mathcal A_n$ before returning to the root. Hence
%$1 -\omega(\rho,\sigma_{1})$ to jump to a child of the root different from $\sigma_{1}$, or to the root. If this happens, then the process hits the root in the next jump. If the process hits $\sigma_{1}$, instead, then $\nu$ is hit before the root with probability
\begin{equation}
\label{eq:q}
q_\omega \df \mathbf{P}_{\omega} ( \widetilde T_{\mathcal A_n}< \widetilde{T}(\rho)) = \sum_{i=1}^b \omega(\rho,\rho^{(i)})\,p_{i,\omega}, \quad \text{where} \quad 
p_{i,\omega}=\Big(\sum_{j=1}^{n} \prod_{k=1}^{r-1} \frac{\omega(\sigma_{k}^{(i)},\sigma_{k-1}^{(i)})}
{\omega(\sigma_{k}^{(i)},\sigma_{k+1}^{(i)})}\Big)^{-1}. 
\end{equation}
 %Define recursively
%$$
%\begin{array}{lll}
% &Z_{1} \df \frac 12 \big(\inf\{n \ge 0 \colon Y_{n} = \sigma_{1}\} + 1\big),  &R_{1} \df \inf\{n \colon Y_{n} = \rho\} ,\\
% &Z_{i} \df   \frac 12 \big(\inf\{n \ge R_{i-1} \colon Y_{n} = \sigma_{1}\} - R_{i-1} + 1\big), &R_i \df \inf\{n >R_{i-1}\colon Y_{n} = \rho\}.     
%\end{array}
%$$
%where by convention $\inf \emptyset=\infty$, and $Z_{i} = \infty$ if $ R_{i-1} = \infty$, for $i\ge 2$. 
%By construction, $Z_i$ counts the number of visits to the root between the successive return times to the root $R_{i-1}$ and $R_i$.  Let 
%$$
%W \df \inf \{ n \ge 1 \colon R_{n} > T(\nu)\}.
%$$
%Notice that 
%$$
%L(\rho, T(\nu)) \le \widetilde{L}(\rho, T(\nu)) \le \sum_{i=1}^{W} Z_{i}\,.
%$$
%and each $Z_{i}$ is independent of $W$ under both the quenched and the annealed measures. In fact they are determined by disjoint parts of the environment. 
It follows that under the quenched measure, $\widetilde L(\rho,\widetilde T_{\mathcal A_n})$ is a geometric variable with parameter $q_\omega$. Hence, with the help of Lemma \ref{lemma:A.2} in the Appendix, we find that
\begin{equation}
 \label{eq:3.1}
 \mathbf E[\widetilde L(\rho,\widetilde T_{\mathcal A_n})^p] \le c_p \mathbb E[q_\omega^{-p}]\,.
\end{equation}
It follows from (\ref{eq:q}), and by independence, that
\begin{equation}
 \label{eq:3.2}
 \mathbb E[q_\omega^{-p}] \le \mathbb E[(\min_i p_{i,\omega})^{-p}\,(1-\omega(\rho,\overleftarrow \rho))^{-p}]
=\mathbb E[(\min_i p_{i,\omega})^{-p}]\,\mathbb E[(1-\omega(\rho,\overleftarrow \rho))^{-p}].
\end{equation}
We use that
\begin{equation}
 \label{eq:3.3}
\mathbb E[(\min_i p_{i,\omega})^{-p}] = \mathbb E[\max_i p_{i,\omega}^{-p}] \le \mathbb E[\Sigma_i p_{i,\omega}^{-p}]=
b\,\mathbb E[p_{1,\omega}^{-p}],
\end{equation}
and we find by (\ref{eq:q}),  by Jensen's inequality and by independence that
\begin{equation}
 \label{eq:3.4}
\mathbb E[p_{1,\omega}^{-p}] \le n^{p-1} \sum_{j=1}^n \mathbb E[A^{-p}]^j= n^{p-1} \frac{\mathbb{E}[A^{-p}]^{n} -1}{\mathbb{E}[A^{-p}] -1}.
\end{equation}
Now observe that 
\begin{equation}
 \label{eq:3.5}
\mathbb E[(1-\omega(\rho,\overleftarrow \rho))^{-p}]=\mathbb E[(1+\tfrac 1{\sum_iA_i})^p],
\end{equation}
and by collecting the results from (\ref{eq:3.1}) to (\ref{eq:3.5}), the claim of the Proposition follows for $n \ge 2$. For $n=1$, a similar (and simpler) argument shows the claim. This finishes the proof of the Proposition.
 \qed
\end{proofsect}

% Let
 
 \begin{proofsect}{Proof of Proposition \ref{prop:momtime}}   
In the course of this proof, we denote with  $\mathbf{Y}^{\ssup \nu}$ the extension of $\mathbf{X}$ to $\overleftarrow{\L}_{\nu}$, and let 
\begin{equation}
\label{dsup}
  D^{\ssup \nu} \df \inf \{ n \ge 1 \colon{Y}^{\ssup \nu}_{n} = \overleftarrow{\nu}\}, \quad \text{and}\quad 
  C(\nu)=\{D^{\ssup \nu}=\infty\}.
\end{equation}
%Set
%$$ 
%C_{n} := \bigcup_{\nu \colon |\nu|=n}  C(\nu).
%$$
Suppose that $|\nu|\ge 1$ and  $C(\nu)$ holds. Then if the process visits $\nu$ it will never return to $\overleftarrow{\nu}$, and in particular it  will not increase the local time spent at the root $\rho$. 
% If $C(\nu)$ holds, after the first time the process hits the first child of $\nu$, if this ever happens,  it will never visit $\nu$ again, and will not increase the local time of $\mathbf Z$  at the root. 
%Roughly, our strategy is to use the extensions on paths to give an upper bound for the total time spent at the root by
% time $T_{k}$ and show that the probability that $\bigcap_{i=1}^{k}C^{c}_{i}$ holds decreases fast enough in $k$.
%Using the independence between disjoint collections of Poisson processes, we infer that  
 Define 
\begin{equation}
\begin{aligned}
d = \inf\{ k \ge 1 \colon &\text{there are $b$ distinct vertices $\nu_1,\ldots,\nu_b$ at level $k$ with different}\\
                          &\text{ancestors at level $1$ s.t.~$C(\nu_i)$ holds for all $1 \le i \le b$}\}.
\end{aligned}
\end{equation} 
On $\{d=n\}$, we choose $b$ distinct vertices $\nu_1,\ldots,\nu_b$ at level $n$ with different ancestors at level $1$ s.t.~$C(\nu_i)$ holds for all $1 \le i \le b$, and in the notation used in Propostion \ref{prop:3.1}, we denote this set of vertices with $\mathcal A_n$. Notice that
\begin{equation}
 \label{eq:3.6}
 L(\rho)\,\1_{\{d=n\}} \le \widetilde L(\rho,\widetilde T_{\mathcal A_n})\,\1_{\{d=n\}}\,.
\end{equation}
%It follows that, in the notation of Proposition \ref{prop:3.1},
%We recall the notation from Proposition \ref{prop:3.1}, and 
With the help of (\ref{eq:3.6}), we infer that for $q,q'$ as in the proposition,
\begin{equation}
 \label{eq:3.7}
 \mathbf E[L(\rho)^p] %=\sum_{n=1}^\infty \mathbf E[L(\rho)^p,d=n] 
 \le \sum_{n=1}^\infty \mathbf E[\widetilde L(\rho,\widetilde T_{\mathcal A_n})^p,d=n]
 \le  \sum_{n=1}^\infty \mathbf E[\widetilde L(\rho,\widetilde T_{\mathcal A_n})^{pq}]^{1/q}\,\mathbf P[d=n]^{1/q'},
\end{equation}
where in the last inequality we used Holder's inequality. Let us now estimate $\mathbf P(d=n)$. 
The events $C(\nu)_{|\nu|=n}$ are determined by disjoint parts of the environment, and are thus independent and identically distributed under the annealed measure. 
%Notice that  $\mathbf{P} (C(\nu)) =1-\beta$, where $\beta$ is defined in \eqref{beta}.
%Observe  that $\{d=1\}$ means that for all vertices $\nu$ at level one the event $C(\nu)$ holds. Hence
%\begin{equation}
% \label{eq:3.8}
% \mathbf P(d=1) = (1-\beta)^b\,.
%\end{equation}
Fix $n \ge 2$. At level $n-1$, there are $b$ families of $b^{n-2}$ vertices each that have different ancestors at level one. If $\{d=n\}$ holds, then the event $C(\cdot)^c$  holds for all 
$b^{n-2}$ vertices in at least one of these families of vertices at level $n-1$. 
With $\mathbf P(C(\cdot))=1-\beta$, it follows
$$
\mathbf P(d=n) \le 1-(1-\mathbf P(C^c)^{b^{n-2}})^b = 1-(1-\beta^{b^{n-2}})^b\,, 
%\le  1-(1-\alpha_\psi^{b^{n-2}})^b
%\le \beta^{b^{n-2}}\,\sum_{i=1}^b (-1)^{i-1} {b \choose i}  \beta^{b^{n-2}(i-1)}\,.
%\{d=n\} \subseteq \{\text{\,within the $b$ collections of $b^{n-2}$ vertices
%there are at most $b-1$ distinct vertices $\nu_i$ at level $n-1$ s.t.~$C(\nu_i)$ holds}\}.
$$
and with Proposition \ref{prop:branching}, it follows that
\begin{equation}
\label{eq:8.5}
\mathbf P(d=n) \le 1-(1-\alpha_\psi^{b^{n-2}})^b
 = \alpha_\psi^{b^{n-2}}\,\sum_{i=1}^b (-1)^{i-1} {b \choose i}  \alpha_\psi^{b^{n-2}(i-1)}\,.
\end{equation}
%Hence
%\begin{equation}
%\label{eq:3.8}
% \begin{aligned}
%\mathbf{P}(d=n) \le& \sum_{i=0}^{b-1} {b^{n-1} \choose i} (1-\beta)^{i} \beta^{b^{n-1}-i} 
%= \beta^{b^{n-1}} \sum_{i=0}^{b-1} {b^{n-1} \choose i} \Big( \frac{1-\beta}\beta \Big)^{i}\,.
%\le& \beta^{b^{n-1}} b^{(b-1)(n-1)} \sum_{i=0}^{b-1}\tfrac 1{i!}\big( \tfrac{1-\beta}\beta \big)^{i}
%\le e^{\tfrac{1-\beta}\beta}\,b^{(n-1)(b-1)}\,\beta^{b^{n-1}}.
%\end{aligned}
% \end{equation}
Together with the trivial bound $\mathbf P(d=1) \le 1$, this finishes the proof of the proposition.
\qed
\end{proofsect}  

%%%%%%%%%%%%%%%%%%%%%%%%%%%%%%%%%%%%%%%%%%%%%%%%%%%%%%%%%%%%%%%%%%%%%%%%%%%%%%%%%%%%%%%%%%%%%%%%%%%%%%%%%%%%%%%%%%%%%%%%%
\subsubsection{An explicit lower bound on the speed and an example}
\noindent Recall $\alpha_\psi$ in Proposition \ref{prop:branching}. 
The propositions \ref{thm:1}, \ref{prop:branching} and  {\ref{prop:momtime} (applied with $p=\varepsilon=1$) imply the following
\begin{theorem}
\label{cor:1}
Assume (\ref{eq:A}), and that $\mathbb E[A^{-2}]<\infty$. Then it holds $\mathbf P$-a.s.~that
\begin{equation*}
\label{eq:lb}
 v \ge \frac{1-\alpha_\psi}{\mathbf E[L(\rho)]} \ge \frac{1-\alpha_\psi}{\theta(2,1)^{1/2}+\sum_{n=2}^\infty \,\theta(2,n)^{1/2}\,\,\alpha_\psi^{b^{n-2}/2}\,\,
\big(\sum_{i=1}^b (-1)^{i-1} {b \choose i} \alpha_\psi^{b^{n-2}(i-1)}\big)^{1/2}}\,>0\,.
\end{equation*}
\end{theorem}
{\bf An example.} Let us now provide an explicit example on the regular binary tree (i.e.~$b=2$). We  choose $A_1=A_2$, and we write $A$ for a copy of $A_1$  resp.~$A_2$. We  choose the following distribution
$$
\mathbb P[A=3/10]=\kappa,\quad \mathbb P[A=7/2]=1-\kappa,\quad \kappa \in (0,1/2].
$$
We compute $m_1$, which is given by the left-hand side of (\ref{eq:5}) with $\psi$ replaced by one, and find that for all $ \kappa \in (0,1/2]$,
\begin{equation}
 m_1=2\,\mathbb E[\tfrac A{1+A}]=(182-128 \kappa)/117>1,
\end{equation}
so that $\alpha_1<1$.
%It follows that $\mathbb E[A^{-1}]=38/21<2$, so in particular we can find $\psi$ such that $m_\psi \ge 1$ (and of course the walk is transient). We now show that it is possible to choose $\psi=1$.
The extinction probability $\alpha_1$ is given by the smallest solution of $x=p_0+p_1x+p_2x^2$, hence $\alpha_1=p_0/p_2$.
Let us compute now the offspring distribution $\mathbf p$. 
We denote the site environment corresponding to the events  $\{A=3/10\}$ and $\{A=7/2\}$ by $\omega_1$ resp.~$\omega_2$. It follows that  
\begin{equation}
\begin{array}{lll}
&\omega_1(\rho,\overleftarrow \rho)=5/8, \quad &\omega_1(\rho, \rho^{(1)})=\omega_1(\rho, \rho^{(2)})=3/16, \\
&\omega_2(\rho,\overleftarrow \rho)= 1/8, \quad &\omega_2(\rho, \rho^{(1)})=\omega_2(\rho, \rho^{(2)})=7/16.
\end{array}
\end{equation}
We obtain that 
\begin{equation}
\begin{aligned}
p_0&=\mathbb E[\omega(\rho,\overleftarrow \rho)]=1/8+\kappa/2,\\
p_1&=\mathbb E[ \omega(\rho, \rho^{(1)})\tfrac{\omega(\rho,\overleftarrow \rho)}{\omega(\rho,\overleftarrow \rho)+\omega(\rho, \rho^{(2)})}]+\mathbb E[ \omega(\rho, \rho^{(2)})\tfrac{\omega(\rho,\overleftarrow \rho)}{\omega(\rho,\overleftarrow \rho)+\omega(\rho, \rho^{(1)})}]
=7/36+11\kappa/117,\\
p_2&=1-p_0-p_1=49/72-139\kappa/234.
\end{aligned}
\end{equation}
%It follows that $m_1=p_1+2p_2=118/117>1$, so that the branching process is supercritical. 
Hence $\alpha_1=(117+468\kappa)/(637-556\kappa)$.
Further we find that $\mathbb E[A^{-2}]-1=4864 \kappa/441-45/49$, and that $\mathbb E[(1+\tfrac 1{2A})^2]=64/49+2560\kappa/441$.
%For $\kappa=1/4$, we find that $\alpha_1=39/83$.
%(for $\kappa=1/2$, we find that $\alpha_1=351/359$). Let us detail the case $\kappa=1/4$.
%From Corollary \ref{cor:1} we find that
%\begin{equation}
%\label{eq:speed}
% v \ge \tfrac 12  (\mathbb E[A^{-2}]-1)^{1/2}
%\mathbb E[(1+\tfrac 1{2A})^2]^{-1/2} e^{\tfrac{1-\alpha_1}{\alpha_1}} \frac {1-\alpha_1}{\sum_{n=1}^\infty (n (\mathbb E[A^{-2}]^n-1)
%2^{n-1} \alpha_1^{2^{n-1}})^{1/2}}\,.
%\end{equation}
Let us choose now $\kappa=1/30$. %We find that $\alpha_1 \approx 0.2144$.
%$\kappa=1/4$ and $\kappa=1/16$. We have that $\alpha_1=39/83$ if $\kappa=1/4$, and 
% $\alpha_1=195/803$ if $\kappa=1/16$.
%We find that $\mathbb E[A^{-2}]=1252/441$, and $\mathbb E[(1+\tfrac 1{2A})^2]=1216/441$. 
%The summands in the infinite sum in Corollary \ref{cor:1} decay double exponentially, and we find that the seventh summand is of the order $10^{-10}$. We therefore only sum up the first seven summands, and find that this finite sum equals $4.66$.
%(when we take the correct $\theta(1,2)$, we obtain $9.21574$). The prefactor equals $0.710918$, so that we obtain
We obtain from Theorem \ref{cor:1} and the above computations that
$$
 v \ge 0.1229 \quad \text{if} \quad \kappa=1/30\,.
$$
%Compare this with the expected drift $\mathbb E[2\omega(\rho,\rho^{(1)})-\omega(\rho,\overleftarrow \rho)]=3/4-\kappa$, which equals $43/60 \approx 0.717$ if $\kappa=1/30$.\\
%{\bf Remark:} For $\kappa=1/16$, we find that $v \ge 0.0934$, and for $\kappa=1/4$, we find $v \ge 0.01317$.\\
%$\alpha_1=195/803$.
%From Proposition \ref{prop:branching} and Theorem \ref{thm:1}, we find that
%\begin{equation}
% v \ge \tfrac 12 e^{\tfrac{1-\alpha_1}{\alpha_1}} (1-\alpha_1) |\mathbb E[A^{-2}]-1|^{1/2}
%\mathbb E[(1+\tfrac 1{2A})^2]^{-1/2} \frac 1{\sum_{n=1}^\infty (n \,|\mathbb E[A^{-2}]^n-1|\,
%2^{n-1} \alpha_1^{2^{n-1}})^{1/2}}\,.
%\end{equation}
%We find that $\mathbb E[A^{-2}]=1252/441$, and $\mathbb E[(1+\tfrac 1{2A})^2]=1216/441$. The summands in the infinite sum in the last display decay double exponentially, and we find that the seventh and the eighth summand is of the order $10^{-8}$ resp.~$10^{-18}$. We therefore only sum up the first seven summands, and find that this finite sum equals $9.48784$ (when we take the correct $\theta(1,2)$, we obtain $9.21574$). The prefactor equals $0.710918$, so that we obtain
%$$
% v \ge 0.710918/9.21574=0.0771417.
%$$
%Compare this with the expected drift $\mathbb E[2\omega(\rho,\rho^{(1)})-\omega(\rho,\overleftarrow \rho)]=1/2$.\\
%{\bf Next try:} Apply H\"older instead of Cauchy-Schwarz in (\ref{eq:3.7})? 
%%%%%%%%%%%%%%%%%%%%%%%%%%%%%%%%%%%%%%%%%%%%%%%%%%%%%%%%%%%%%%%%%%%%%%%%%%%%%%%%%%%%%%%%%%%%%%%%%%%%%%%%%%%%
\subsection{Once edge-reinforced random walk}

Durrett, Kesten and Limic \cite{DKL2002} prove transience and provide a law of large numbers with positive speed for once edge-reinforced random walk on a regular tree. However their methods do not give a lower bound for the speed that is always positive. Collevecchio \cite{C2006b} proves  transience for this process defined on supercritical Galton--Watson trees. The same was proved, independently and with different methods by Dai \cite{D2005}. In this section, we provide a lower bound on the speed by using a refinement of the methods from \cite{C2006b}.\\  
%We now provide a lower bound for the speed of once edge-reinforced random walks on regular trees. 
Let us first define the process. Fix $\delta >0$, and denote with $\{\nu,\mu\}$ the edge connecting the neighboring vertices $\nu$  and $\mu$. 
Once $\delta$-edge-reinforced random walk (ORRW($\delta$) or simply ORRW) $\mathbf{X}= \{ X_k, k \ge 0 \}$ is a discrete-time process 
%with values on the vertices of 
on the regular $b$-ary tree $\Gcal_{b}$, and is defined as follows. Each edge  has initial weight one,  i.e.~$W(\{\nu,\mu\},0)=1$, 
with the exception of the edge $\{\overleftarrow{\rho},\rho\}$, which has weight $\delta$,
i.e.~$W(\{\overleftarrow{\rho},\rho\},0)=\delta$. 
This exception helps to simplify our exposition. This initial weight configuration is called {\it initially fair}. 
For $n \ge 1$, we update the weight $W$ of the edges according to the following rule:
\begin{equation}
\label{eq:5.0}
W(\{\nu,\mu\},n)=
\begin{cases}
\delta, \,&\text{ if $\{X_{k-1},X_k\}=\{\nu,\mu\}$ for some $1 \le k \le n$},\\
%\1_{\{T(\nu) \le n \text{ and }T(\mu) \le n\}}+ W(\{\nu,\mu\}, 0)\1_{\{T(\nu) > n \text{ or }T(\mu) > n\}},
1,\, &\text{ otherwise.}
\end{cases}
\end{equation}
%where the stopping times $T(\cdot)$ are defined as in \eqref{T}.
%for $k\ge 0$ and some vertex $\nu$, that $X_{k}=\nu$. 
%
%\noindent Define 
%\[ z_{k}(\nu, \overleftarrow{\nu})= z_k(\overleftarrow{\nu},\nu) \df \left\{ \begin{array}{ll}
%         \delta  & \mbox{  if  $ T(\nu)\le k  $  }\\
%                               \\
%         1 & \mbox{  otherwise } \end{array} \right. . \]
%\{ X_{j},  1 \le j \le k\}$, 
%let $G_{n}(\eta)= \{ X_{j} = \eta $ for some $ j \le n\}$. We introduce 
%$z_{\eta,n}=\delta\, \1_{G_{n}(\eta)} +\1_{G^{c}_{n}(\eta)}$, and for a neighbor $\mu$ of $X_n$, we define
ORRW starts from  $\rho$, i.e.~$X_0 = \rho$, and we define inductively $\mathcal F_n=\sigma(X_0,X_1,\ldots,X_n)$, and the transition probabilities 
\begin{equation}
\label{eq:5.0a}
\mathbf{P}( X_{n+1}=\mu \mid \mathcal F_n ) = \frac{W(\{X_n,\mu\},n)}{\sum_{\nu: \nu \sim X_n} W(\{X_n,\nu\},n)},
\end{equation} 
if $\mu$ is a neighbor of $X_n$, and zero otherwise. The canonical law of this process is denoted with $\mathbf P$.
Later on, we will also use the following initial weights, where not only the edge $\{\overleftarrow{\rho},\rho\}$ has weight $\delta$, but a connected collection of edges containing the edge $\{\overleftarrow{\rho},\rho\}$, i.e.~if some edge has weight $\delta$, then each edge on the path connecting this edge to the root has weight $\delta$.
%So far, we considered  initial conditions
% $W(\{\nu, \mu\}, 0\}) =1$ for $\{\nu, \mu\} \neq\{\overleftarrow{\rho}, \rho\}$ and  equal to $\delta$ otherwise. Next we define ORRW with different initial conditions. 
We denote with $\mathbb{W}$ the set of such initial weight configurations.
Of course, $\mathbb W$ contains the initially fair weights, that we denote from now on with $w_0$.
%connected sets of edges containing the edge $\{\overleftarrow{\rho}, \rho\}$. We denote with $w_{0}$  the set containing only the edge $\{\overleftarrow{\rho}, \rho\}$.  
For $w \in \mathbb W$ let $w(\{\nu, \mu\})$ be the weight  that $w$ assigns to the edge $\{\nu,\mu\}$.
For any weight configuration $w \in  \mathbb{W}$, define  $W_{w}(\{\nu,\mu\},0) = w(\{\nu,\mu\})$, and for $n \ge 1$,
%\delta \1_{\{ \{\nu,\mu\} \in w\}} + \1_{\{ \{\nu,\mu\} \notin w\}},$ and 
$$
W_{w}(\{\nu,\mu\},n)=
\begin{cases}
\delta, \,&\text{ if $\{X_{k-1},X_k\}=\{\nu,\mu\}$ for some $1 \le k \le n$},\\
%\1_{\{T(\nu) \le n \text{ and }T(\mu) \le n\}}+ W(\{\nu,\mu\}, 0)\1_{\{T(\nu) > n \text{ or }T(\mu) > n\}},
w(\{\nu,\mu\}),\, &\text{ otherwise.}
\end{cases}
$$ 
%Define  $\mathbf{X}^{\ssup w} $ to be  ORRW  with initial configuration $w$,  i.e. a process which starts at $\rho$ and with transition probability $\mathbf{P}_{w}$ given by
The transition probabilities are defined similarly as in \eqref{eq:5.0a}, with $W(\cdot,n)$ replaced by $W_w(\cdot,n)$.
%$$ 
%\mathbf{P}_{w}( X^{\ssup w}_{n+1}=\mu \mid \mathcal F_n ) = \frac{W_{w}(\{X^{\ssup w}_n,\mu\},n)}{\sum_{\nu: \nu \sim X^{\ssup w}_n} W_{w}(\{X^{\ssup w}_n,\nu\},n)}.
%$$
The canonical law of ORRW started at $\rho$ and in the initial weight configuration $w \in \mathbb W$ is denoted with $ \mathbf{P}_w$ (clearly $\mathbf{P}=  \mathbf{P}_{ w_{0}})$. 
%and $\mathbf{X}= \mathbf{X}^{\ssup {w_{0}}}$.
 Recall the exponential random variables $h_{k}(\cdot,\cdot),\,\, k \ge 1$,  with mean one, used in definition~\ref{ext} and fix a subtree $\Ccal$ of $\Gcal_b$.

\begin{definition}\rm[Extension $\mathbf{Y}^{\Ccal}$ on the subtree $\Ccal$]
\label{def:3}
The extension $\mathbf{Y}^{\Ccal}$ of $\mathbf{X}$ on the subtree $\Ccal$ is defined as follows. 
Fix a starting point $\eta$ in $\Ccal$, i.e.~$Y^\Ccal_0=\eta$ and an initial weight configuration $w\in \mathbb{W}$. We define $\mathbf Y^\Ccal$ iteratively in the following way.
 Let $ s_{1}(\nu) $ be the first time $\mathbf{Y}^{\Ccal}$
  reaches some vertex $\nu$. Define $ N_\nu^\Ccal$ to be the set of neighbors of $\nu$ in $\Ccal$.
 The first jump after $s_{1}(\nu)$  is towards the neighbor $\mu \in N_{\nu}^\Ccal$ for which the following minimum
\begin{equation} 
\label{expjump1}
\min_{\mu \in N_\nu^\Ccal} \frac{h_{1}(\nu,\mu)}{W_{w}(\{\nu,\mu\},s_1(\nu))}
%h_{1}(\nu,\overleftarrow{\nu}) \wedge \min_{i}  \big(h_{1}(\nu,\overrightarrow{\nu}^{\ssup i})/A^{\ssup i}_{\nu}\big),
\end{equation}
%where $a \wedge b= \min\{a,b\}$. 
is a.s.~attained.
We define  $s_{k}(\nu),\,\, k \ge 2$, inductively via
\begin{eqnarray*}
 && s_{k}(\nu) \df  \inf \big\{ n > s_{k-1} \colon Y^{\Ccal}_{n} = \nu \big\}, \mbox{ and} \\
 &&  j_{k}(\nu,\mu)  \df 1+ \mbox{ number of times $\mathbf{Y}^{\Ccal}$ jumped from $\nu$ to its  neighbor $\mu$  by time $s_{k}$}.
 \end{eqnarray*}
The first jump after $s_{k}(\nu)$  is towards the neighbor $\mu$ for which the following minimum 
\begin{equation}
\label{expjump2}
\min_{\mu \in N_{\nu}^\Ccal} \frac{h_{j_{k}}(\nu,\mu)}{W_{w}(\{\nu,\mu\},s_1(\nu))}
%h_{j(\nu,\overleftarrow{\nu})}(\nu,\overleftarrow{\nu})\wedge \min_{i} \big( h_{j(\nu, \overrightarrow{\nu}^{\ssup i})}(\nu,\overrightarrow{\nu}^{\ssup i})/A^{\ssup i}_{\nu} \big)
\end{equation}
is a.s.~attained.
\qed \end{definition}
%\begin{remark}\rm
The comments in remark \ref{extcoupl} also apply here.
%Notice that  for any subtree $\Ccal$, the jumps of the  extension $\mathbf{Y}^{\Ccal}$ and those  of $\mathbf{X}$ on $\Ccal$ are determined via the same collection of exponentials $h_{i}(\cdot, \cdot)$. \qed
%\end{remark}
\noindent We now introduce a similar color scheme as in definition \ref{def:color}.
\begin{definition}\rm
\label{def:color2}
  Fix an integer $\psi \ge 1$, and denote with 
$\mathbf Y(\overleftarrow \mu,\nu)$, for a descendant $\nu$ of $\mu$, the extension of ORRW on the ray connecting $\overleftarrow \mu$ to $\nu$, started at $\mu$, in the following initial weight configuration. The edge $\{\overleftarrow \mu, \mu\}$ has weight $\delta$ and all the other edges in the path connecting $\mu$ to $\nu$ have initial weight 1. 
%Now we apply the same color scheme as in definition \ref{def:color}.
%\qed
A vertex $\nu$ at level $\psi$ is colored if and only if  $\mathbf Y(\overleftarrow \rho,\nu)$
hits $\nu$ before  $\overleftarrow \rho$.  A vertex $\nu$ at level $k\psi,\,k \ge 2$, is colored if and only if
\begin{itemize}
\item its ancestor at level $(k-1)\psi$, say $\mu$, is colored, and 
\item  $Y(\overleftarrow \mu,\nu)$ hits $\nu$ before  $\overleftarrow \mu$. 
\end{itemize}
All the other vertices are uncolored, and only vertices that are at a level $k\psi$, $k\ge 1$, can be colored. 
\qed \end{definition}
%Since $\mathbf Y(\overleftarrow \rho,\nu)$ is started at $\overleftarrow \rho$, the first jump is towards $\rho$, so that $W(\{\overleftarrow \rho,\rho\},1)=\delta$. We could as well start $\mathbf Y(\overleftarrow \rho,\nu)$ at $\rho$, and assign initial weight $\delta$ to the edge $\{\overleftarrow \rho,\rho\}$.
\noindent This color scheme constitutes again a homogeneous branching process, with extinction probability $\alpha_\psi$. 
Notice that for every $b \ge 2$, and every $\delta>0$, we can always find an integer $\psi \ge 1$ such that 
\begin{equation}
\label{eq:5.1}
b^{\psi} \prod_{j=1}^{\psi} \frac j{j + \delta}>1.
\end{equation}
We define $D$ in the same way as in \eqref{beta0}, and also $\beta_w=\mathbf P_w(D=\infty)$, and we write $\beta=\beta_{w_0}$. Recall $\mathbb W$ below \eqref{eq:5.0a}. We have the following
\begin{proposition}
\label{lemma:weights}
If $\psi$ is such that \eqref{eq:5.1} holds, then $\alpha_\psi <1$.  
If $\delta>1$, then for every $w \in \mathbb W$, it holds that 
%For all $w \in \mathbb W$, it holds that 
$\beta_w \le \alpha_\psi$.
% On the other hand, if $\delta \le 1$, then for every $w \in \mathbb W$, it holds that 
%%For all $w \in \mathbb W$, it holds that 
%$\beta_w \le \alpha_\psi(1).$
%If the colored process associated to $\mathbf X$ survives, then the colored process associated to $\mathbf Z$ survives.
\end{proposition}
\begin{proofsect}{Proof}
%Consider the weight configuration $w_0=\delta^{\otimes 1} \times 1^{\otimes \psi+1}$ on the ray connecting $\overleftarrow \rho$ to a vertex $\nu$ at level $\psi$. 
The probability that  $\mathbf Y(\overleftarrow \rho, \nu)$, started at $\rho$, in the initially fair weight configuration $w_0$, 
hits level $\psi$ before it hits $\overleftarrow \rho$ is equal to (see Lemma 1 in \cite{C2006b})
\begin{equation}%\label{ruinor}
    \prod_{j=1}^{\psi} \frac j{j + \delta}.
\end{equation}
%and equality holds if the initial weights assigned to each edge are all equal 1.
Hence the mean of the offspring distribution of the colored process is equal to $b^{\psi} \prod_{j=1}^{\psi} \frac j{j + \delta}$, which is larger than one by our choice of $\psi$. This shows that $\alpha_\psi<1$.
%If the colored process survives, then $D=\infty$. It follows that $\beta \le \alpha_\psi$.
%If the colored process survives, then $D=\infty$. Hence $\mathbf P(D<\infty) \le \alpha_\psi$.
 Now choose an  initial weight configuration $w \in \mathbb W$. 
%and denote with $\mathbf Z$ the ORRW associated to $w$ and started at $\rho$.
% Since the initial weights in $w$ for $\mathbf Z$ are dominating the initial weights for $\mathbf X$,
If $\delta >1$, we can couple the extension  $\mathbf Y(\overleftarrow \rho, \nu)$, started at $\rho$, in the initially fair weight configuration $w_0$, to the extension 
$\widetilde{\mathbf Y}(\overleftarrow \rho, \nu)$, started at $\rho$, in the weight configuration $w \in \mathbb W$,  in such a way that 
$|\widetilde{\mathbf Y}| \ge  |\mathbf Y|$.  %To see this, notice that  both of these processes start from $\rho$. 
To do this, we choose a family of independent variables $(E_n^{\uparrow},E_n^\downarrow)_{n \ge 1}$, with i.i.d.~exponential entries with mean 1. At each time point $n$, the vector $(E_n^{\uparrow},E_n^\downarrow)$ is attached {\it both} to the positions $Y_n$ and $\widetilde{Y}_n$, with $E_n^{\uparrow}$ attached to the edge connecting $Y_n$ and $\widetilde{Y}_n$ to the vertex $\nu$ at level $|Y_n|+1$ resp.~$\widetilde \nu$ at level $|\widetilde{Y}_n|+1$, and $E_n^{\downarrow}$ attached to the edge connecting $Y_n$ and $\widetilde{Y}_n$ to the vertex $\mu$ at level $|Y_n|-1$ resp.~$\widetilde \mu$ at level $|\widetilde{ Y}_n|-1$. The jump of $\mathbf Y$ at time $n+1$ is to the vertex $\nu$ or $\mu$ for which the minimum
\begin{equation}
\label{expjump3}
\min \{\,\,\frac{E_n^\uparrow}{W_{w_0}(\{Y_n,\nu\},n)}, \frac{E_n^\downarrow}{W_{w_0}(\{Y_n,\mu\},n)}\,\,\}
%h_{j(\nu,\overleftarrow{\nu})}(\nu,\overleftarrow{\nu})\wedge \min_{i} \big( h_{j(\nu, \overrightarrow{\nu}^{\ssup i})}(\nu,\overrightarrow{\nu}^{\ssup i})/A^{\ssup i}_{\nu} \big)
\end{equation}
is a.s.~attained, and similarly for $\widetilde{\mathbf Y}$, where we replace the weights $W_{w_0}$ by $W_w$, and the vertices $\nu$, $\mu$ by $\widetilde{\nu}$, $\widetilde{\mu}$. Notice that in this way the extensions $\mathbf Y$ and $\widetilde{\mathbf Y}$ have the same distribution as in the definition \ref{def:3}.
Let 
$$
r=\inf\{n \ge 1 \colon |Y_n| \neq |\widetilde{Y}_{n}|\}
$$
 be the first splitting time, and for ease of notation, 
 let  $e_{0}, e_{1}$ be the two edges incident to $Y_{r-1}=\widetilde{Y}_{r-1}$, where $e_{1}$ connects $Y_{r-1}$ to its child on the path, and $e_{0}$ connects $Y_{r-1}$ to its parent $\overleftarrow{Y}_{r-1}$.
Clearly $W_w(e_0,r-1)=W_{w_0}(e_0,r-1)=\delta$, since the edge $e_0$ is crossed by both processes. Also, by construction, $W_{w}(e,r-1) \ge W_{w_{0}}(e, r-1)$ for any edge $e$ lying on the path connecting $\overleftarrow \rho$ to $\nu$. If we would have $W_{w}(e_{1},r-1) = W_{w_{0}}(e_{1}, r-1)$, then, by the construction of the coupling in \eqref{expjump3}, $Y_{r}=\widetilde{Y}_{r}$, a contradiction. Hence $W_{w}(e_{1},r-1)=\delta$  and $W_{w_{0}}(e_{1}, r-1)=1$. It follows again from \eqref{expjump3} that the only way $\mathbf Y$ and $\widetilde{\mathbf Y}$ can split is that $|\widetilde{Y_{r}}|= |Y_{r}| +2$. Define 
$$
s=\inf\{n > r \colon |Y_n|=|\widetilde{Y}_{n}|\}.
$$
For any edge $e$ lying on the path connecting $\overleftarrow \rho$ to $\nu$,  we have that  $W_{w}(e,s) \ge W_{w_{0}}(e,s)$, and we can reiterate the previous argument to prove  that $|\widetilde{\mathbf Y}| \ge  |\mathbf Y|$.
%so that $Y_{r-1} = \widetilde{Y}_{r-1}$ but  $Y_{r} \neq \widetilde{Y}_{r}$. Let  $e_{0}, e_{1}$ be the two edges incident to $Y_{r-1}$, where $e_{1}$ connects $Y_{r-1}$ to the ancestor of $\nu$ which lies at level $|Y_{r-1}|+1$, and $e_{0}$ connects $Y_{r-1}$ to  $\overleftarrow Y_{r-1}$.  Recall that $W_{w}(e_{1},r-1) \ge W_{w_{0}}(e_{1}, r-1)$, and $ W_{w}(e_{1},r-1), W_{w_{0}}(e_{1}, r-1)\in \{1,\delta\}$. As $Y_{r-1} = \widetilde{Y}_{r-1}$, we have $W_{w}(e_{0},r-1) = W_{w_{0}}(e_{0}, r-1)= \delta$. 
%As the two processes are generated using the same exponentials, if we had $W_{w}(e_{1},r-1) = W_{w_{0}}(e_{1}, r-1)$, then the two process would jump to the same vertex at time $r$,  contradicting our definition of $r$.  Hence $ W_{w}(e_{1},r-1) > W_{w_{0}}(e_{1}, r-1)$, which implies that  soon after time $r-1$, the process $\mathbf{Y}$ jumps across $e_{0}$ while $\widetilde{\bf Y}$ jumps through $e_{1}$. Hence  $|\widetilde{Y_{r}}|= |Y_{r}| +2$. Let $t$ be the first time the two processes meet again.  For any edge $e$ lying on the path connecting $\overleftarrow \rho$ to $\nu$,  we have that  $W_{w}(e, t) \ge W_{w_{0}}(e,t)$, and we can reiterate the previous  argument to prove  that $|\widetilde{\mathbf Y}| \ge  |\mathbf Y|$.
Consider the coloring process, defined in the same way as above \eqref{eq:5.1}, but on the weight configuration $w$. It follows that, if the  coloring process associated to $\mathbf Y$ survives,  then as $|\widetilde{\mathbf Y}| \ge  |\mathbf Y|$, the coloring process  
associated to  $\widetilde{\mathbf Y}$ survives. But on this last event, $D=\infty$. 
Hence  $\beta_w=\mathbf P_w(D<\infty) \le \alpha_\psi$. 
%In the following, we denote ORRW(1) with $\mathbf Z$. Of course $\mathbf Z$ is simple random walk. 
%If $\delta \le 1$, 
%we can couple the extension  $\mathbf Z(\overleftarrow \rho, \nu)$ to $\widetilde{\mathbf Y}(\overleftarrow \rho, \nu)$ in the  weight configuration $w$ such that $|\widetilde{\mathbf Y}| \ge |\mathbf Z|$. 
%Now the coloring process associated to $\mathbf Z$ survives if for instance $\psi \ge 2$, 
%so that the coloring process associated to $\mathbf Y$ in the weight configuration $w$ survives. As above, it follows that 
%$\beta_w \le \alpha_\psi(1)$.
%This finishes the proof.
%This means that, if the colored process associated to $\mathbf X$ survives, then the colored process associated to $\mathbf Z$ survives. In other words, 
%$\mathbf P_{w}(D<\infty) \le \mathbf P(D<\infty)$, and the claim follows.
\qed
\end{proofsect}
The random variable $L(\cdot)$ is defined in the same way as in \eqref{eq:L}. We have the following
\begin{proposition} 
\label{prop:L}
If $\delta>1$, under $\mathbf P_{w_{0}}$, the random variable $L(\rho)$ is stochastically dominated by  a geometric variable with parameter $(1-\alpha_\psi)\,b/(b+\delta)$.
%\begin{equation}
%(1-\alpha_\psi)\,\frac b{b+\delta}.
%\end{equation}
\end{proposition}
\begin{proofsect}{Proof} %Let $\mathbf{Y}^{\ssup \rho}$ be the extension on $\overleftarrow{\L_{\nu}}$, started at $\overleftarrow{\nu}$.  
Recall that $\mathbf{X}$ starts from $\rho$ in the initially fair weight configuration $w_0$.
%and that $\{\overleftarrow \delta,\delta\}$ is the only edge with weight  $\delta$, and all the other edges have initial weight equal to 1. 
With probability $b/(b+\delta)$ the first  jump will be towards one of the children of $\rho$. Then, started at this child of $\rho$, with probability $1-\beta$, the process will never return to $\rho$. Whenever it returns  to $\rho$,  it starts  on some random weight configuration $w \in \mathbb W$, depending on the past of the path. 
Under $\mathbf P_w$, the probability that ORRW jumps to one of the children of $\rho$ is greater than  $b/(b+\delta)$.
To see this, recall that the edge $\{\overleftarrow \rho,\rho\}$ has weight $\delta$, and we change all the weights on the edges connecting $\rho$ to its children to one. Since $\delta >1$, this decreases the probability to jump to level one, and we obtain the lower bound for this probability. 
% since some of the edges connecting $\rho$ to its children have weight  $\delta>1$, which increases the probability to jump to level one.
 Under $\mathbf P_w$, ORRW, started at a child $\nu$ of $\rho$, 
 has probability larger than $1-\beta_{\overline w}$ of never returning  to $\rho$, where $\overline{w}$ is the weight configuration induced by $w$ on $\overleftarrow{\L}_{\nu}$. 
%More precisely $\{\rho, \nu\} \in \overline{w} $ and for any other edge $e$ connecting two vertices in  $\overleftarrow{\L}_{\nu}$, we have that  $e \in \overline{w}$ if and only if $ e \in w$.  
With the help of Proposition \ref{lemma:weights}, we find that, for any $w \in \mathbb W$, the escape probability from $\rho$ is at least 
$(1-\alpha_\psi)b/(b+\delta)$, and it follows that the number of returns to $\rho$ is stochastically dominated  by a geometric variable with parameter $ (1-\alpha_\psi)b/(b+\delta)$. 
%If $\delta \le 1$, under $\mathbf P_w$, the probability that ORRW jumps to one of the children of $\rho$ is greater than  $b/(b+1)$. As above, the edge $\{\overleftarrow \rho,\rho\}$ has weight $\delta$, and we change all the weights on the edges connecting $\rho$ to its children to $\delta$. Since $\delta <1$, this decreases the probability to jump to level one. Proceeding as above, we find by Proposition \ref{lemma:weights} that the number of returns to $\rho$ is stochastically dominated  by a geometric variable with parameter $ (1-\alpha_\psi(1))\,b/(b+1)$. 
%The proof is finished.
\qed
\end{proofsect}
We recall from \cite{DKL2002} that a law of large numbers with positive speed holds, i.e.~$\mathbf P$-a.s., $v = \lim_{n \to \infty}|X_{n}|/n>0$. Further it is shown that $v \le b/(b+\delta)$, but no lower bound is available.
We are now ready to provide a lower bound for the speed that is always positive.  
\begin{theorem}
\label{thm:ORRW}
If $\delta >1$, choose $\psi \ge 1$ such that \eqref{eq:5.1} holds. Then the speed $v$ satisfies
\begin{equation}
\label{eq:6.1}
v \ge \frac{1-\beta}{\mathbf E[L(\rho)]} \ge (1-\alpha_\psi)^{2}\frac b{b+\delta}>0.
\end{equation}
\end{theorem}
\begin{remark}\rm
\label{rem:A.1a}
%Actually Theorem \ref{thm:ORRW} holds for $\delta >1$. 
Notice that in the case of $\delta<b$ we can compare $|\mathbf{X}|$ with a simple random walk on the non-negative integers with drift equal to $(b- \delta)/(b+\delta)>0$. It follows that for $\delta <b$ we have $v \ge (b-\delta)/(b+\delta)$. In this case, we find that the lower bound in \eqref{eq:6.1} is larger than $(b-\delta)/(b+\delta)$ if and only if $\alpha_\psi < 1-\sqrt{1-\delta/b}$. The challenging case is $\delta \ge b$, which is covered by Theorem \ref{thm:ORRW}.
% It might be that $\alpha_\psi(1)=1$, see also the example below for $b=2,\psi=1$. There $\alpha_1(1)=0$ since the condition for survival of the colored process is $\delta<1$. If we choose $\psi=3, b=3$, the condition is $\delta<2$, hence there $\alpha_1(1)<1$. The best thing would be to have a method that gives the same lower bound in the case $\delta \le 1$ as in the case $\delta >1$.
\end{remark}
\begin{proofsect}{Proof of Theorem \ref{thm:ORRW}}
%The same result as Lemma \ref{lemma:A.1} in the previous section holds, with exactly the same proof. 
Define the random variable $\Pi_{k}$ in the same way as in \eqref{eq:Pi}. Observe that the same result as Lemma \ref{lemma:A.1} in the previous section holds, with exactly the same proof. By straightforward modifications, we further see that Proposition \ref{thm:1} holds in the setting of once edge-reinforced random walk. The first inequality follows.
%by a straightforward modification as in the proof of Proposition \ref{thm:1}. 
The second and third inequality then follow directly from Propositions \ref{prop:L} and \ref{lemma:weights}.
%\begin{Lemma}
%\label{lemma:pi2}
%The random variable $\Pi_{k}$ is stochastically dominated by a  geometric random variable  with parameter $\mathbf P(D=\infty)$.
%\end{Lemma}
%The proof is identical to the proof of Lemma \ref{lemma:A.1} in the previous section.
%The claim now follows from Proposition \ref{prop:L} and Lemma \ref{lemma:pi2}.
\qed
\end{proofsect}
Next we show monotonicity of the lower bound on the speed in \eqref{eq:6.1}.
\begin{proposition}%[Monotonicity of the lower bound]
Choose $\delta_2>\delta_1 \ge 1$. Then for every $\psi \ge 1$, $\alpha_\psi(\delta_1) \le \alpha_\psi(\delta_2)$, and in particular the lower bound in \eqref{eq:6.1} is decreasing in $\delta$ for $\delta >1$.
\end{proposition}
\begin{proofsect}{Proof}
Denote with $\mathbf Y^{(1)}$ and $\mathbf Y^{(2)}$ the extensions on rays $[\overleftarrow \rho,\infty)$ corresponding to ORRW($\delta_1$) resp.~ORRW($\delta_2$), started at $\rho$, in the initially fair weight configuration $w_0^{(\delta_1)}$ resp.~$w_0^{(\delta_2)}$. 
%Clearly $W_{w_0^{(\delta_2)}}(e,n) \ge W_{w_0^{(\delta_1)}}(e,n)$ for all edges $e$ and all times $n$. 
Using the same coupling as in \eqref{expjump3}, we can show that 
%Using the same exponentials to run the two processes,  we have  
$|\mathbf Y^{(1)}| \ge |\mathbf Y^{(2)}|$. To see this, call $r$ to be the first time the two processes split, and let $e_{0}$ and  $e_{1}$ be as in in proof of Proposition ~~\ref{lemma:weights}. Next we show that none of the processes traversed edge $e_{1}$ by time $r-1$.   In fact, as the two processes coincide up to time $r-1$, if one of them traversed $e_{1}$, also the other did. On the other hand, both of them traversed $e_{0}$ by time $r-1$, in order to reach $Y_{r-1}^{(1)}= Y_{r-1}^{(2)}$.  Hence $\mathbf{P}(|Y_{r}^{(1)}|= |Y_{r-1}^{(1)}|+1) = 1/2 = \mathbf{P}(|Y_{r}^{(2)}|= |Y_{r-1}^{(2)}|+1)$. By construction of the coupling, this would imply that $Y_{r}^{(1)} = Y_{r}^{(2)}$, which contradicts the definition of $r$.
As none of the processes traversed edge $e_{1}$ by time $r-1$, while both traversed $e_{0}$, using the fact  $\delta_{2}>\delta_{1}$ we  infer that $|Y_{r}^{(1)}| > |Y_{r}^{(2)}|$.
  Denote with $t$ the first time, after $r$, when  the two processes meet, and
%We have that $\mathbf{P}(|Y_{t+1}^{(1)}|= |Y_{t}^{(1)}|+1) = 1/2 \ge \mathbf{P}(|Y_{t+1}^{(2)}|= |Y_{t}^{(2)}|+1)$. 
let $r_{1}$ be the first time after $t$, when the two processes split again.
As $|Y_{k}^{(1)}| \ge |Y_{k}^{(2)}|$ for all $ k \le r_{1}-1$, we have that there is no edge reinforced by $Y_{k}^{(2)}$  which has not been reinforced by $ Y_{k}^{(1)}$, $k \le r_{1}-1$. 
This, together with the fact that $\delta_1 > \delta_{2}>1$, and the construction of the coupling, implies that  $\mathbf{P}(|Y_{r_{1}}^{(1)}|= |Y_{r_{1}-1}^{(1)}|+1)  \ge\mathbf{P}(|Y_{r_{1}}^{(2)}|= |Y_{r_{1}-1}^{(2)}|+1)$. By construction of the coupling, we have that $ |Y_{r_{1}}^{(1)}| > |Y_{r_{1}}^{(2)}|$. By reiterating this argument, we get $|\mathbf Y^{(1)}| \ge |\mathbf Y^{(2)}|$.
This implies that for every $\psi$, $\alpha_\psi(\delta_1) \le \alpha_\psi(\delta_2)$, and it follows that the lower bound in \eqref{eq:6.1} is decreasing in $\delta$.
\qed
\end{proofsect}
\section{Moment bounds on the first regeneration time}

In addition to providing an explicit lower bound on the speed, our methods can be extended to give an explicit upper bound on the tail of a certain regeneration level. We present a unified approach that applies both for random walk in a random environment and once edge-reinforced random walk. Hence, in what follows, $\mathbf X$ denotes either one of these processes. We start by defining the regeneration times.
\begin{definition}\rm
\label{def:1} 
We define the first regeneration level as  follows
%the first time the walk hits a level for the first time and then never backtracks to the previous level:
$$
\ell_1 \df\inf \{k \ge 1\colon  D(X_{T_k})=\infty\},
$$
and iteratively
$$
\ell_n \df \inf \{k > \ell_{n-1}\colon  D(X_{T_k})=\infty   \},
$$
where $D(\cdot)$ is defined in (\ref{beta0}) and we use the convention $\inf \emptyset = \infty$. The regeneration times  are defined as $\tau_n=T_{\ell_{n}}, n \ge 1$, on the event  $\{\ell_n <\infty\}$. 

 %Level $j \ge 1 $ is a  regeneration level if  the first jump after $T_{j}$ is towards level $j+1$, and after time $T_{j+1}$  the process never backtracks to the parent of $X_{T_j}$.   
 %Define $\ell_1$ to be the  regeneration level with minimum distance from the root, and for $i >1$,
%\[ \ell_i := \min \{ j > \ell_{i-1}\colon j \mbox{ is a regeneration level}\}. \]
%\item Define the $i$-th regeneration time to be $\tau_i  \df T_{\ell_i} $. Notice that $\ell_i = | X_{\tau_i}|$.
\qed \end{definition}
%\begin{definition}
%\label{def:1} 
%We define the first regeneration time as the first time the walk hits a level for the first time and then never backtracks to the previous level:
%$$
%\tau_1 \df \inf \{k \ge 1: k=T_{|X_{T_k}|}, D(X_{T_k})=\infty\},
%$$
%and iteratively
%$$
%\tau_n \df \inf \{k \ge \tau_{n-1}: k=T_{|X_{T_k}|}, D(X_{T_k})=\infty \},
%$$
%with the convention $\inf \emptyset = \infty$. The regeneration levels are defined as $\ell_n=|X_{\tau_n}|, n \ge 1$, on the event  $\{\tau_n <\infty\}$.
 %Level $j \ge 1 $ is a  regeneration level if  the first jump after $T_{j}$ is towards level $j+1$, and after time $T_{j+1}$  the process never backtracks to the parent of $X_{T_j}$.   
 %Define $\ell_1$ to be the  regeneration level with minimum distance from the root, and for $i >1$,
%\[ \ell_i := \min \{ j > \ell_{i-1}\colon j \mbox{ is a regeneration level}\}. \]
%\item Define the $i$-th regeneration time to be $\tau_i  \df T_{\ell_i} $. Notice that $\ell_i = | X_{\tau_i}|$.
%\qed \end{definition}
In other words, a regeneration time occurs   when the walk hits a level for the first time and then never backtracks to the previous level. Clearly, these are not stopping times.
It is easy to see that under transience, it holds that for all $n \ge 1$, $\tau_n <\infty$ $\mathbf P$-a.s.
%, see \cite{G2004}. 
It is also known that in the setting of random walks in random environment, the first regeneration level $\ell_1$ has exponential moments under the conditioned measure $\mathbf P(\cdot|D=\infty)$. This is for instance proved in in Lemma 4.2 in \cite{DGPZ2002} for biased random walks on Galton-Watson trees, and can be directly adapted to our setting. 
For once edge-reinforced random walk with $\delta>1$, we know that $\ell_1$ has all moments finite under  $\mathbf P(\cdot|D=\infty)$, see Lemma 7 in \cite{DKL2002} (this statement is actually proved for certain cut levels, but notice that our regeneration level is smaller than the cut level in \cite{DKL2002}).\\
We now present a unified approach that applies to both settings, and that provides explicit estimates for the tail of $\ell_{1}$ and for the moments of $\tau_{1}$.

\subsection{The tail of the first regeneration level}
We assume that 
we can choose $\psi$ such that \eqref{eq:A} is fulfilled for random walk in
 a random environment resp.~once edge-reinforced random walk.  
Recall that for ORRW($\delta$), this  is always possible, see (\ref{eq:5.1}) and Proposition \ref{lemma:weights}.\\
We will find {\it explicit} exponential tails on $\ell_1$.  These tail estimates on $\ell_1$ are obtained by refining the color scheme from definitions \ref{def:color} resp.~\ref{def:color2}.
\begin{definition}\rm
Let $\nu$ be a vertex at level $k\psi,\,k \ge 1$. 
%For any vertex $\nu $ such that $|\nu|$ is a multiple of $\psi$ we now define certain sets of vertices $\Theta_{\nu}$
%and $\Sigma_{\nu}$ that are subsets of the tree $\Lambda_\nu$.  
Let $\Theta_{\nu}$ be the set of vertices $\mu$ in $\L_{\nu}$ which are first children and whose distance from $\nu$   is a multiple of $\zeta \psi$.
%Let $\Theta_{\nu}$ be the set of vertices $\mu$ such that
%\begin{itemize}
%\item $\mu$ is a descendant of $\nu$, 
%\item the difference $|\mu|-|\nu|$ is a multiple of $\zeta \psi$,
%\mathbb N=\{\zeta \psi,2\zeta \psi,\ldots\}$
%\item $\mu$ is a first child,
%\end{itemize}
%%and their descendants.
Let  $\Sigma_{\nu}$ be the set of vertices $\mu$ in  $\Lambda_{\nu}$ %(the subtree attached to $\nu$) 
such that
\begin{itemize}
%\item $\mu$ is a descendant of $\nu$, 
\item $\mu$ is colored (in particular $|\mu|$ - $|\nu|$ is a multiple of $\psi$),
\item all ancestors of $\mu$ in $\Lambda_{\nu}$ do not belong to $\Theta_{\nu}$. 
\end{itemize}
Further define
% \begin{equation}
% \label{sub1}
% \begin{aligned}
$ B(\nu)= \{ \Sigma_{\nu} \mbox{ is infinite}\},  \text{ and }  B_0 = B(\rho),\,\,\, 
 B_{i} = B(X_{T_{\psi \zeta i}}), \,\,i \ge 1$.
% \end{aligned}
% \end{equation}
\qed \end{definition}
In other words, $\Sigma_\nu$ is the set of colored vertices in $\Lambda_\nu$ minus the colored vertices that are elements of subtrees generated by vertices $\mu$ that are first children and  $|\mu|-|\nu|=k\zeta \psi,\,k \ge 1$. 
In a first step, we introduce an auxiliary branching process and use it to derive an explicit lower bound on the probability of $B_0$, see Lemma \ref{return1}.  In a second step, in Lemma \ref{aind},  we then show that the events $B_i$ are independent. 
In \cite{C2009}, section 3, the counterpart of these lemmata for vertex-reinforced jump processes are stated and proved in a similar way.\\
%Denote with $\tilde T_i$ the (thinned out) sequence of regeneration times such that $|X_{\tilde T_i}| \in \zeta n \mathbb N$ for all $i \ge 1$. By convention $\tilde T_0=T_0=0$.
For any pair of distributions $f_1$ and $f_2$, denote by $f_1 \,\overline{*}\,f_2$ the distribution of $\sum^{V}_{ k = 1} M_{k}$, where 
\begin{itemize}
\item $V$ has distribution  $f_1$, and 
\item $\{ M_{k},\, k \in  \N \}$ is a sequence of i.i.d.~random variables, independent of $V$, each with distribution $f_2$.
\end{itemize}
We set $\mathbf{p}^{\ssup 1} :=  \mathbf{p}$, and define, by recursion, $\mathbf{p}^{\ssup j} :=  \mathbf{p}^{\ssup {j-1}} \,\overline{*}\, \mathbf{p}$ for $j \ge 2$. The distribution $\mathbf{p}^{\ssup j}$ describes the number of elements, at time $j$, in a population which evolves like a branching process generated by one ancestor and  with offspring distribution $\mathbf{p}$. 
%If we let
%\[ 
%m \df \sum_{j=0}^{2^n} j p_{j},
%\]
%The mean of $\mathbf{p}^{\ssup j}$, for $j \ge 2$,  is $m^{j}$, where $m$ was introduced in (\ref{eq:m}). Recall that $m>1$.\\
Let $q_{0} = p_{0} + p_{1},\, $ and for $k \in \{1,\ldots,b^\psi-1\}$, set    $ q_{k} = p_{k+1}$. Set $\mathbf{q}$ to be the distribution which assigns to $i \in \{0, \ldots, b^\psi-1\}$ probability $q_{i}$.  For $j \ge 2$, let
$\mathbf{q}^{\ssup j} := \mathbf{p}^{\ssup {j-1}} \,\overline{*}\, \mathbf{q}$. Denote by $q^{\ssup j}_{i}$ the weight that the distribution $ \mathbf{q}^{\ssup j}$ assigns to $i \in \{ 0, \ldots, (b^\psi-1)b^{(j-1)\psi}\}$. The mean of $\mathbf{q}^{\ssup j}$ is $m_\psi^{j-1}(m_\psi-1)$. From now on, $\zeta$ denotes the smallest positive integer    such that 
\begin{equation}\label{offwthfc}
m_\psi^{\zeta-1}(m_\psi-1) > 1.
\end{equation}
(This is possible since we chose $\psi$ such that $m_\psi >1$.)
Define $\gamma $ to be the smallest positive solution of the equation
\begin{equation}
\label{root}
 x = \sum_{k=0}^{\vartheta} x^{k}q^{\ssup \zeta}_{k}, \quad \text{where $ \vartheta=b^{(\zeta-1)\psi}(b^\psi-1)$}.
\end{equation}
%while  $\zeta$ and $(q_{k}^{\ssup \psi})$ have been defined at the beginning of this section.
\begin{Lemma}
\label{return1}
 Assume \eqref{offwthfc}. We have that for $i \ge 0$,
 %\begin{equation}\label{ret1}
 $\mathbf{P}(B_{i})=\mathbf{P}(B_{0}) \ge 1-\gamma>0$. 
 %%\qquad \forall i \in n \N.
 %\end{equation}
\end{Lemma}

 \begin{proofsect}{Proof}
Fix $i$ and notice that by stationarity, $\mathbf P(B_i)=\mathbf P(B_0)$.
 From the definition of $\Sigma_\rho$, it follows that the offspring distribution of colored vertices at level $\zeta \psi$ in  $\Sigma_\rho$ is obtained as follows. The number of vertices at level $(\zeta-1)\psi$ 
has law $\mathbf{p}^{\ssup{(\zeta-1)\psi}}$. Each vertex at level  $\zeta \psi$ has a number of colored offspring distributed as 
 $\mathbf{p}=\mathbf{p}^{\ssup 1}$. If from each of these offspring we  delete the first child, the number of the remaining colored offspring is distributed as $ \mathbf{q}$. 
Hence the offspring distribution modeling $\Sigma_\nu$ is given by $\mathbf{q}^{\ssup {\zeta}}=\mathbf{p}^{\ssup {\zeta-1}}\,\overline{*}\, \mathbf{q}$.
Then, from the basic theory of branching processes we know that 
the 
extinction probability equals the smallest positive solution of the equation \eqref{root}. In virtue of \eqref{offwthfc}  we have that $\gamma<1$.
 \qed           \end{proofsect}

\begin{Lemma}
\label{aind}
 The events $B_{i}$, $i \ge 1$,  are independent under $\mathbf P$.
 \end{Lemma}
 \begin{proofsect}{Proof}
%We remark that $\zeta \ge 2$. 
Choose  integers $ 0 < i_{1}<i_{2}< \ldots<i_{k}$. % with $i_{j} \in \zeta n \N := \{ \zeta n, 2 \zeta n, 3 \zeta n, \ldots\}$ for all $j \in \{1, 2, \ldots, k\}$.  
It is enough to prove that
\begin{equation}
\label{indep1}
\mathbf{P} ( \bigcap_{j=1}^{k} B_{i_{j}} ) = \prod_{j=1}^{k} \mathbf{P} ( B_{i_{j}}).
\end{equation} 
We proceed by backward recursion.
%Fix a vertex $ \nu$ at level $i_{k}$. 
We use the notation introduced in definition~~\ref{ext}.
The set $ B(\nu)$ belongs to the sigma-algebra generated by $\{h_{i}(\eta, \mu) \colon \; \eta,\mu \in \mbox{Vert}(\L_{\nu}) $ and $ i \ge 1 \}$.
%$\{  P(u,w) \colon \; u,w \in \mbox{Vert}(\L_{\nu}) \}$. 
Notice that each $X_{T_i},\,i \ge 1$, is a first child. Hence the set 
$\cap_{j=1}^{k-1} B_{i_{j}} \cap \{X_{T_{\psi\zeta i_{k}}} = \nu\}$ belongs to $\{  h_{i}(\eta, \mu)  \colon \; \eta \notin \mbox{Vert}(\L_{\nu}) \}$.
%$\{  P(u,w) \colon \; u \notin \mbox{Vert}(\L_{\nu}) \}$. 
As the two events belong to disjoint collections of independent exponential variables, they are independent. 
We have
\begin{equation*}
\label{eq:2.10}
\begin{aligned}
 \mathbf{P} ( \bigcap_{j=1}^{k} B_{i_{j}} )
%\mathbf{P} \big( B_{i_{k}} \cap \bigcap_{j=1}^{k-1} B_{i_{j}} \big) 
= \sum_{\nu} \mathbf{P}\big(B_{i_{k}}  \cap  \bigcap_{j=1}^{k-1} B_{i_{j}} \cap \{X_{T_{\psi\zeta i_{k}}} = \nu\}\big)
%&= \sum_{\nu} \mathbf{P}\Big(B(\nu)  \cap  \bigcap_{j=1}^{k-1} B_{i_{j}} \cap \{X_{T_{i_{k}}} = \nu\}\Big)
= \sum_{\nu} \mathbf{P}(B(\nu))\,\,\mathbf{P} \big(  \bigcap_{j=1}^{k-1} B_{i_{j}} \cap \{X_{T_{\psi\zeta i_{k}}} = \nu\} \big).  
%&= \mathbf{P}\big(B_0 \big) \sum_{\nu} \mathbf{P} \Big(  \bigcap_{j=1}^{k-1} B_{i_{j}} \cap \{X_{T_{i_{k}}} = \nu\} \Big)
%=&\mathbf{P}(B_0)\,\, \mathbf{P} \big(  \bigcap_{j=1}^{k-1} B_{i_{j}} \big).
\end{aligned}
\end{equation*}
From  stationarity, it follows that $\mathbf P(B(\nu))=\mathbf P(B_0)$, and from the independence of $B(\nu)$ and $\{X_{T_{i\psi\zeta}} = \nu\}$, we infer that for an arbitrary vertex $\nu$, and each $i \ge 1$,
%and in virtue of the self-similarity property of the regular tree 
\begin{equation}
\mathbf P(B(\nu))=\mathbf{P}(B_{i})\,.
%\sum_\nu  \mathbf{P}\big(B(\nu),\,X_{T_{i_k}}=\nu \big)
%=\sum_\nu \mathbf{P}\big(B_0\big)\,\mathbf P \big(X_{T_{i_k}}=\nu \big)=\mathbf P\big(B_0\big).
\end{equation}
Now the right-hand side of (\ref{eq:2.10}) equals 
\begin{equation}
\mathbf P(B_0) \sum_{\nu} \mathbf{P} \big(  \bigcap_{j=1}^{k-1} B_{i_{j}} \cap \{X_{T_{\psi\zeta i_{k}}} = \nu \}\big)=\mathbf{P}(B_{i_k})  \mathbf{P} \big(  \bigcap_{j=1}^{k-1} B_{i_{j}} \big)\,.
%\sum_\nu  \mathbf{P}\big(B(\nu),\,X_{T_{i_k}}=\nu \big)
%=\sum_\nu \mathbf{P}\big(B_0\big)\,\mathbf P \big(X_{T_{i_k}}=\nu \big)=\mathbf P\big(B_0\big).
\end{equation}
(\ref{indep1}) follows now by iteration.
% Hence
%\begin{equation}\label{recstep}
%\mathbf{P} \Big( B_{i_{k}} \cap \bigcap_{j=1}^{k-1} B_{i_{j}} \Big)= \mathbf{P}\big(B_{i_{k}} \big) \mathbf{P} \Big(  \bigcap_{j=1}^{k-1} B_{i_{j}} \Big).
%\end{equation}
%Reiterating \eqref{recstep} we get \eqref{indep1}.
\qed
  \end{proofsect}

\begin{theorem}
\label{codaelle}
 Assume (\ref{eq:A}). For $n \ge 1$, we have that
\begin{equation}\label{tailelle1}
 \mathbf{P} ( \ell_{1}  \ge n \psi \zeta  ) \le \gamma^{n-1}, 
\end{equation}
where $\gamma$ is defined in (\ref{root}). 
\end{theorem}
\begin{proofsect}{Proof}  
Notice that on the event $B_{i}$, the colored  process survives in the subtree $\Lambda_{X_{T_{i \psi \zeta}}}$. It follows that  $B_i \subseteq \{\text{level $i \psi \zeta$ is a regeneration level}\}$. Hence 
$$
\{ \ell_{1}  \ge n \psi \zeta  \} \subseteq \bigcap_{i=1}^{n-1} B_i^c\,,
$$
and the Theorem now follows from the Lemmata \ref{return1} and \ref{aind}.
%A level which is of type ${ \bf A}$ is a regeneration level. 
\qed           \end{proofsect} 

%%%%%%%%%%%%%%%%%%%%%%%%%%%%%%%%%%%%%%%%%%%%%%%%%%%%%%%%%%%%%%%%%%%%%%%%%%%%%%%%%%%%%%%%%%%%%%%%%%%%%%%%%%%%%%%%%%%
\subsection{Moment bounds for the first regeneration time}

Recall the first regeneration time in Definition \ref{def:1}, and define
$$
\Pi =\sum_{\nu \in \Lambda} \1_{\{T(\nu) \le \tau_1\}}
$$
%\sum_{n=1}^{\infty}\sum_{k=1}^{n}\Pi_{n,k}\1_{\{\ell_1=n\}},
to be the number of distinct vertices visited by time $\tau_{1}$.
%Further let $\nu_{k,i},i=1,\ldots,b^k,$ be an enumeration of the vertices at level $k$, and let 
%$$
%\Pi_{k}=\sum_{i=1}^{b^k}\1_{\{T_{\nu_{k,i}}<\infty\}}
%$$ 
%be the number of vertices visited at level $k$.
We denote with $M(n,q)$ the $n$-th moment of a geometric variable with parameter $q$. 
We have the following explicit bound on the moments of $\Pi$, which implies an explicit bound on the moments of $\tau_1$, see Theorem \ref{thm:tau} below. 
\begin{proposition}
\label{prop:Pi}
Assume (\ref{eq:A}). For $p \ge 1$, it holds that
\begin{equation}
\label{eq:8.0}
\mathbf{E}[\Pi^p] \le \gamma^{-1/2}\Big(1-\gamma^{\tfrac 1{2\psi\zeta}}\Big)^{-1}\,\Big(M\big(p,1-\gamma^{\tfrac 1{2\psi\zeta}}\big)-1\Big)\,M^{1/2}(2p,1-\beta) \,.
%\le c_p c_{2p}^{1/2} (1-\beta)^{-p} \,\,\gamma^{-1/2} \Big(1-\gamma^{\tfrac1{2\psi\zeta}}\Big)^{-p-1}.
%c_{2p}^{1/2}(1-\beta)^{-p}\sum_{n=1}^{\infty}(n+1)^{p} \,\gamma^{\tfrac n{2\psi\zeta}-\tfrac 12}\;<\infty,
\end{equation}
\end{proposition}
\begin{proofsect}{Proof}
Recall $\Pi_k=\sum_{\nu \colon |\nu|=k} \1_{\{T(\nu) < \infty\}},\, k \ge -1$, which is the number of vertices visited at level $k$, and observe that
\begin{equation}
 \Pi = \sum_{n=1}^\infty \sum_{\nu} \1_{\{T(\nu) \le T_n\}}\,\1_{\{\ell_1=n\}} 
\le \sum_{n=1}^\infty \sum_{\nu: |\nu| <n} \1_{\{T(\nu) < \infty\}}\,\1_{\{\ell_1=n\}}
= \sum_{n=1}^\infty \sum_{k=-1}^{n-1} \Pi_k\,\1_{\{\ell_1=n\}}.
\end{equation}
We use Jensen's inequality, and obtain that
\begin{equation}
\label{eq:8.1}
 \mathbf{E}[\Pi^p] \le \sum_{n=1}^{\infty} \mathbf{E}[\big(\sum_{k=-1}^{n-1}\Pi_{k}\big)^p \1_{\{\ell_1=n\}}]
 \stackrel{(Jensen)}{\le} \sum_{n=1}^{\infty}(n+1)^{p-1}\sum_{k=-1}^{n-1}\mathbf{E}[\Pi_{k}^p \1_{\{\ell_1=n\}}].
\end{equation}
First notice that Lemma \ref{lemma:A.1}, proved for random walk in a random environment, holds also for once edge-reinforced random walk with the same proof. We first use Cauchy-Schwarz's inequality, and then Lemma \ref{lemma:A.1} together with Lemma \ref{lemma:A.2} from the Appendix to obtain that the right-hand side of the last display is smaller than
\begin{equation}
\label{eq:8.2}
 \sum_{n=1}^{\infty}(n+1)^{p-1}\sum_{k=-1}^{n-1}\mathbf{E}[\Pi_{k}^{2p}]^{1/2}\,\mathbf{P}( \ell_1=n)^{1/2}
 \le %c_{2p}^{1/2}(1-\beta)^{-p}
 M^{1/2}(2p,1-\beta) \sum_{n=1}^{\infty}(n+1)^{p}\,\mathbf{P}(\ell_1 \ge n)^{1/2}.
\end{equation}
Finally, with Theorem \ref{codaelle}, we obtain that
\begin{equation}
\label{eq:8.3}
 \sum_{n=1}^{\infty}(n+1)^{p}\,\mathbf{P}(\ell_1 \ge n)^{1/2} \le 
 \gamma^{-\tfrac 12} \sum_{n=2}^\infty n^p \,\gamma^{\tfrac {n-1}{2\psi\zeta}}
=\gamma^{-1/2}\Big(1-\gamma^{\tfrac 1{2\psi\zeta}}\Big)^{-1}\,\Big(M(p,1-\gamma^{\tfrac 1{2\psi\zeta}})-1\Big).
\end{equation}
The claim \eqref{eq:8.0} now follows by collecting the results in  \eqref{eq:8.1} to \eqref{eq:8.3}.
%claim of the proposition follows.
\qed 
\end{proofsect}
\noindent We are now ready to state the main result of this subsection.
\begin{theorem} 
\label{thm:tau}
Assume (\ref{eq:A}) and that $\mathbb{E}[A^{-p-\varepsilon}]<\infty$ for some $p \ge 1$ and $\varepsilon >0$.  %$\mathbb{E}[\omega(\rho, \overrightarrow{\rho}^{\ssup 1})^{-4p}]<\infty$, 
It holds that
\begin{equation}
\mathbf{E}[\tau_{1}^{p}] \le \tfrac{\pi^{2}}{6}\,\mathbf{E}[L(\rho)^{p+\varepsilon}]^{\tfrac 1q}\,
                             \mathbf{E}[ \Pi^{2(p-1)q'} ]^{\tfrac 1{2q'}}\,\mathbf{E}[\Pi^{4q'}]^{\tfrac 1{2q'}}<\infty,
%\mathbf{E}[\tau_{1}^{p}\;|\; D=\infty] \le \frac 1{1-\alpha} \frac{\pi^{2}}{6}\Theta(4p)^{1/2} \Xi(4(p-1))^{1/4}\Xi(8)^{1/4},
\end{equation}
where $q=1+\varepsilon/p$, and $q'=1+p/\varepsilon$ is the dual of $q$.
%where we recall 
%\begin{equation}
%\begin{aligned}
%\Theta(p)=&\sum_{\heap{F \in \Pcal_{b}}{F \neq \emptyset}}  c(p) \mathbf{E}\Big[\big(\sum_{\nu \in F}\omega(\rho,\nu) \big)^{-p}\Big] \beta^{b - |F|} (1- \beta)^{|F|}\\ &+  \big(c(p) \mathbb{E}[\omega^{-p}(\rho, \overrightarrow{\rho}^{\ssup 1})] \big)^{1/2} \Big(\sum_{k=2}^{\infty}  \Big(1 + c(p)k^{p-1}\frac{\mathbb{E}[A^{-p}]^{k} -1}{\mathbb{E}[A^{-p}] -1}\Big)^{1/2}   \beta^{b^{k-1}/2}\Big),\\
%\mbox{and}&\\
% \Xi(p) =& c(2p)^{1/2}(1-\beta)^{-p}\sum_{n=1}^{\infty}(n+1)^{p} \,
%\gamma^{\tfrac n{2\psi\zeta}-\tfrac 12}.
%\end{aligned}
%\end{equation}
\end{theorem}
%\begin{remark}
%\label{rem:tau}
%\rm
%Taking into account remark \ref{rem:momtime}, and applying again H\"older's inequality with adjusted parameters  in \eqref{eq:4.2}, we find that we can relax 
%\end{remark}
\begin{proofsect}{Proof}
By Jensen's inequality, we find
\begin{equation}
\label{eq:4.1}
\begin{aligned}
\mathbf{E}[\tau_{1}^{p}] &= \mathbf{E} [ \big(\sum_{i=1}^{\Pi} L(\sigma_{i})\big)^{p}]
\le \mathbf{E}[ \Pi^{p-1} \sum_{i=1}^{\Pi} L(\sigma_{i})^{p}] 
= \sum_{i=1}^{\infty}  \mathbf{E}[ \Pi^{p-1} L(\sigma_{i})^{p} \1_{\{\Pi \ge i\}}].
\end{aligned}
\end{equation}
By H\"older's inequality, and by stationarity, the right-hand side of the last display is smaller than
\begin{equation}
\label{eq:4.2}
\begin{aligned}
 %&\sum_{i=1}^{\infty} \mathbf{E}[L(\sigma_{i})^{2p}]^{1/2} \mathbf{E}[ \Pi^{4(p-1)}]^{1/4} \mathbf{P}(\Pi \ge i)^{1/4} \\
&  \mathbf{E}[L(\rho)^{p+\varepsilon}]^{\tfrac 1q}\,\mathbf{E}[ \Pi^{2(p-1)q'} ]^{\tfrac 1{2q'}} \,\sum_{i=1}^{\infty} \mathbf{P}(\Pi \ge i)^{\tfrac 1{2q'}}
%\le& \mathbf{E}[L(\rho)^{2p}]^{1/2}\, \mathbf{E}[ \Pi^{4(p-1)}]^{1/4}\,\sum_{i=1}^{\infty}i^{-2}\mathbf{E}[\Pi^{8}]\\
%=&\tfrac{\pi^{2}}{6} \mathbf{E}[L(\rho)^{2p}]^{1/2}\, \mathbf{E}[ \Pi^{4(p-1)}]^{1/4} \mathbf{E}[\Pi^{8}]^{1/4}.
\end{aligned}
\end{equation}
By Chebychev's inequality, we find that
\begin{equation}
\label{eq:4.3}
\sum_{i=1}^{\infty} \mathbf{P}(\Pi \ge i)^{\tfrac 1{2q'}} \le \sum_{i=1}^{\infty}i^{-2}\mathbf{E}[\Pi^{4q'}]^{\tfrac 1{2q'}}
=\tfrac{\pi^{2}}{6}\mathbf{E}[\Pi^{4q'}]^{\tfrac 1{2q'}}.
\end{equation}
Putting (\ref{eq:4.1}),(\ref{eq:4.2}) and (\ref{eq:4.3}) together, we obtain the claim.
%The result follows by using Propositions  \ref{momtime} and \ref{prop:Pi}, and the fact that
%$$ \mathbf{E}[\tau_{1}^{p}\;|\; D= \infty] \le \frac{\mathbf{E}[\tau_{1}^{p}]}{\mathbf{P}(D=\infty)}\le \frac{\mathbf{E}[\tau_{1}^{p}]}{1- \alpha} .$$
\qed
\end{proofsect}
%Proposition \ref{prop:momtime} provides an explicit upper bound on $\mathbf E[L(\rho)^p]$. To obtain an explicit upper bound on $\mathbf E[\tau_1^p]$, we need the following  result.
%%%%%%%%%%%%%%%%%%%%%%%%%%%%%%%%%%%%%%%%%%%%%%%%%%%%%%%%%%%%%%%%%%%%%%%%%%%%%%%%%%%%%%%%%%%%%%%%%%%%%%%%%%%%%%%%%%%%%%%%
\subsection{An invariance principle and bounds on the covariance}

For ORRW, an invariance principle is known, see Theorem 3 in Durrett, Kesten and Limic \cite{DKL2002}. For RWRE, an annealed invariance principle easily follows from the results of Aid\'ekon \cite{A2009}. We further refer to Peres and Zeitouni \cite{PZ2008} for a quenched invaraince principle for biased random walks on Galton-Watson trees.
%We denote with $D[0,\infty)$ the set of real-valued functions on $[0,\infty)$, which are right-continuous and have left limits (also known as c\`adl\`ag functions). We endow this set with the Skorohod topology and its Borel-$\sigma$-algebra, cf.~\cite{ek} p.116. We introduce the $D[0,\infty)$-valued random variable
%\begin{equation}
%B^n_\cdot=\frac 1{\sqrt n}(|X_{[\cdot n]}-[\cdot n]v),\,\,\,n \ge 1,
%\end{equation}
%where $[t]$ denotes the integer part of $t$.
Define 
$$
B^n_\cdot=\frac 1{\sqrt n}(|X_{[\cdot n]}|-[\cdot n]v),\,\,\,\beta^n_t=B^n_t+(nt-[nt])(B^n_{t+1}-B^n_t),\,\,\,n \ge 1,
$$  
i.e.~$\beta$ 
%and let  the $C(\mathbb R_+,\mathbb R)$-valued random variable
%\beta^n_t=B^n_t+(nt-[nt])(B^n_{t+1}-B^n_t),
is the polygonal interpolation of $ k/n \to B^n_{ k/n},\,\,k \ge 0$. 
We endow the space $C(\mathbb R_+,\mathbb R)$ of continuous functions with the topology of uniform convergence on compacts, and with its Borel $\sigma$-algebra. 
\begin{proposition}
\label{thm:3}
%If (\ref{eq:A}) holds and if $\mathbb{E}[A^{-2-\varepsilon}]<\infty$ for some $\varepsilon >0$, then 
The %$D[0,\infty)$-valued random variable $B^n_\cdot$ and the 
$C(\mathbb R_+,\mathbb R)$-valued random variable $\beta^n_\cdot$ converge under $\mathbf P$ in law to a Brownian motion $B_\cdot$ with covariance 
$$ 
K= \mathbf{E}[(\ell_{1}- v \tau_{1})^{2}|D=\infty]\,\, \mathbf{E}[\tau_{1}|D=\infty]^{-1}.
$$
%In addition we have the following  bounds on $K$:
%\begin{equation}
% \label{eq:4.4}
% K 
%\end{equation}
\end{proposition}
\begin{proofsect}{Proof}
%Following Proposition 10.4 p.149 in \cite{ek}, the convergence in law of $B^n\cdot$ and $\beta_\cdot^n$ are equivalent.
%We  show that $\beta^n_\cdot$ converges in law under $\mathbf P$ to a Brownian motion with covariance matrix $K$. 
For ORRW, we refer to Theorem 3 in \cite{DKL2002}. For RWRE, observe that the second moment of $\tau_1$, und thus of $\ell_1$, is finite, as follows from Propositions 2.1 and 2.2 in Aid\'ekon \cite{A2009}. Since $\mathbf P[D=\infty]=1-\beta >0$, also $\mathbf E[\ell^2_1|D=\infty] \le \mathbf E[\tau^2_1|D=\infty]<\infty$. Further it is well-known that  
\begin{equation}
\label{eq:iid}
\begin{aligned}
&(\tau_{i+1}- \tau_{i}, \ell_{i+1}-\ell_{i})_{i \ge 1} \text{ is an i.i.d.~sequence under $\mathbf{P}$, and for $i\ge 1$,}\\
&\text{$(\tau_{i+1}- \tau_{i}, \ell_{i+1}-\ell_{i})$ has same law under $\mathbf P$ as $(\tau_{1}, \ell_{1})$ under $\mathbf{P}(\,\cdot\,|\, D= \infty)$},
\end{aligned}
\end{equation}
 see \cite{G2004} (see also \cite{LPP1996} for a similar statement for biased random walks on Galton-Watson trees).
With the help of this i.i.d.~structure, the proof of the invariance principle is now quite standard, see for instance Theorem 3 in Durrett, Kesten and Limic \cite{DKL2002} and also Theorem 3.3 in Shen \cite{Shen2003}. 
\qed
\end{proofsect}
With the help of Theorem \ref{codaelle} and Theorem \ref{thm:tau}, we obtain explicit bounds on the covariance $K$ via the following proposition. 
For RWRE (resp. ORRW)  denote with  $w$ the right-hand side in inequality \eqref{eq:lb} (resp. (\ref{eq:6.1})) , so that $ v \ge w$. Let  $a$   be the smallest even integer larger or equal to $ [3/w]+1$.  As $w \le 1$, we have  $a \ge 4$.
\begin{proposition}
\label{prop:5}
In the case of RWRE, we assume that 
 (\ref{eq:A})  holds and that $\mathbb{E}[A^{-2-\varepsilon}]<\infty$ for some $\varepsilon >0$.  In the case of ORRW we choose $\psi$ satisfying (\ref{eq:5.1}).
%If (\ref{eq:A}) resp.~(\ref{eq:5.1}) hold and if $\mathbb{E}[A^{-2-\varepsilon}]<\infty$ for some $\varepsilon >0$, 
Then we have the following common upper bound on the covariance $K$
\begin{equation}
\label{eq:4.57}
\begin{aligned}
K \le \mathbf (1-\alpha_\psi)^{-1}(\mathbf E[\ell_1^2]+\mathbf E[\tau_1^2])\quad \text{for RWRE and ORRW},
\end{aligned}
\end{equation}
and the following lower bound %for random walk in random environment and once edge-reinforced random walk respectively
\begin{equation}
\label{eq:4.57a}
\begin{aligned}
K &\ge b\,(1-\alpha_\psi)\,\mathbf E[\tau_1]^{-1}\,\mathbb{E}[\omega(\rho,\, \overrightarrow{\rho}_{1})^{\frac a2}] \mathbb{E}[\omega( \overrightarrow{\rho}_{1},\, \rho)^{\frac a2-1} \big(1- \omega( \overrightarrow{\rho}_{1},\, \rho)\big)] \quad \text{for RWRE},\\
K &\ge (1-\alpha_\psi)\,\mathbf E[\tau_1]^{-1}\,\Big(\frac b{b+\delta}\Big)^2\, \Big(\frac \delta{b+\delta}\Big)^{a/2-1}\,\Big(\frac \delta{b-1+2\delta}\Big)^{a/2-1}\quad \text{for ORRW}.
\end{aligned}
\end{equation}
\end{proposition}
\begin{proofsect}{Proof of Proposition \ref{prop:5}}
We start with the upper bound. We use the trivial bound $(a - b)^{2} \le a^{2} + b^{2},\, a, b \ge 0$,  and $v \le 1$ to obtain that
\begin{equation}
 \label{eq:4.58}
 K \le  \mathbf E[\ell_1^2|D=\infty] + \mathbf E[\tau_1^2|D=\infty] \le (1-\beta)^{-1} (\mathbf E[\ell_1^2] + \mathbf E[\tau_1^2]).
\end{equation}
The upper bound  \eqref{eq:4.57} follows from Proposition \ref{prop:branching}. Let us now turn to the lower bound \eqref{eq:4.57a} for random walk in random environment. We use the following approach
\begin{equation}
 \label{eq:4.59}
%\begin{aligned}
%K&= \mathbf{E}[\tau_{1}|D=\infty]^{-1}\\
\mathbf{E}[(\ell_{1}- v \tau_{1})^{2}|D=\infty] 
\ge   \mathbf{E}[(\ell_{1}- v \tau_{1})^{2} \1_{\{ v \tau_{1} \ge \ell_{1} +1 \}}|D=\infty] 
\ge \mathbf{P}[ v \tau_{1} \ge \ell_{1} +1 |D=\infty],
%\end{aligned}
\end{equation}
where the last inequality comes from the fact that on the event $\{v \tau_{1} \ge \ell_{1} +1\}$ we have  $ (\ell_{1}- v \tau_{1})^{2} \ge 1$. Hence
\begin{equation}
 \label{eq:4.60}
%\begin{aligned}
K \ge \mathbf{P}( v \tau_{1} \ge \ell_{1} +1| D=\infty)\,\, \mathbf{E}[\tau_{1}|D=\infty]^{-1}
\ge \mathbf{P}( v \tau_{1} \ge \ell_{1} +1, D=\infty)  \mathbf{E}[\tau_{1}]^{-1}.
%\end{aligned}
\end{equation}
Next we find a suitable subset of $ \{v \tau_{1} \ge \ell_{1} +1\} $ whose probability is easy to compute. 
Consider the event  
$$ C \df \{ T_{2} = a,\,D(X_{T_2})=\infty,\,\cup_{i=1}^b\{X_{j} \in \{\rho, \overrightarrow{\rho}_{i}\}, \, \forall j \le T_{2}-1\}\}.
$$
If this event holds then the walk, started at the root $\rho$, visits level two first at time $a$ and, after this time,  never goes back to level 1.Moreover before time $T_{2}$, the process $\mathbf{X}$ visits only the vertices $\rho $ and $ \overrightarrow{\rho_{i}}$ for some $i$, and
hence it does not return to $\overleftarrow{\rho}$.  
As $a \ge 4$, it jumps at least once from $\overrightarrow{\rho_{1}}$ to $\rho$, so that level one cannot be a cut level and $\ell_{1} =2$.   As $ a \ge [3/w]+1 \ge [3/v]+1$, we have 
\begin{equation*}
 C \subset \{\ell_{1} =2, \tau_{1} \ge [3/v]+1,\,D=\infty\}.
\end{equation*}
On the event $\{\ell_{1} = 2, \tau_{1} \ge  [3/v]+1 \}$ we have that  $v \tau_{1} \ge 3$, hence $v \tau_{1} - \ell_{1} \ge 1$. In other words, 
\begin{equation}
\label{eq:4.60a}
C \subset \{v \tau_{1} \ge \ell_{1} +1,\,D=\infty\}.
\end{equation}
We first focus on the RWRE case. Let us now compute the probability of the event $C$. The Markov property implies that
%Given the event $\{ T_{2} = a\} \cap\{X_{j} \in \{\rho, \overrightarrow{\rho}_{1}\}, \, \forall j \le T_{2}-1\}$, the process $\mathbf{X}$ makes his first jump towards $ \overrightarrow{\rho}_{1}$.   The process keeps oscillating between the two vertices $ \rho$ and  $ \overrightarrow{\rho}_{1}$ until time $a$ when it first hits a child of $  \overrightarrow{\rho}_{1}$. Under the quenched measure, the probability of the event  $\{ T_{2} = a\} \cap\{X_{j} \in \{\rho, \overrightarrow{\rho}_{1}\}, \, \forall j \le T_{2}-1\}$ equals to 
$$
\mathbf P_\omega(C)=\sum_{i=1}^b \omega(\rho, \overrightarrow{\rho_{i}})^{\frac a2} \omega( \overrightarrow{\rho_{i}},\rho)^{\frac a2-1} \big(1- \omega( \overrightarrow{\rho_{i}},\, \rho)\big)\,\mathbf E_\omega[\mathbf P_{X_{T_2},\omega}(D=\infty)].
$$
The random variables $\omega(\rho, \overrightarrow{\rho_{i}})$,  $\omega( \overrightarrow{\rho_{i}}, \rho)\big(1- \omega( \overrightarrow{\rho_{i}},\, \rho)\big)$ and  $\mathbf E_\omega[\mathbf P_{X_{T_2},\omega}(D=\infty)]$ are independent, since they are measurable w.r.t.~
 disjoint parts of the environment. We use in addition stationarity to find that
$$
\mathbf P(C)=b\, \mathbb E[\omega(\rho, \overrightarrow{\rho_1})^{\frac a2}]\,\mathbb E[ \omega( \overrightarrow{\rho_1}, \rho)^{\frac a2-1}(1- \omega( \overrightarrow{\rho_1}, \rho))]\,\mathbf E[\mathbf P_{X_{T_2},\omega}(D=\infty)].
$$
Again, by independence and stationarity,
$$
\begin{aligned}
\mathbf E[P_{X_{T_2},\omega}(D=\infty)] &=\sum_\nu \mathbf E[\mathbf P_{\nu,\omega}(D=\infty),\,X_{T_2}=\nu]\\&=
\sum_\nu \mathbf P_{\nu}(D=\infty)\,\mathbf P(X_{T_2}=\nu)=\mathbf P(D=\infty)=1-\beta.
\end{aligned}
$$
It follows that
\begin{equation}
\label{eq:4.61}
\mathbf{P}(v \tau_{1} - \ell_{1} \ge 1, D= \infty) \ge  \mathbf{P}(C) =b\, \mathbb{E}[\omega(\rho, \overrightarrow{\rho_1})^{\frac a2}] \mathbb{E}[\omega( \overrightarrow{\rho_1}, \rho)^{\frac a2-1} \big(1- \omega( \overrightarrow{\rho_1},\rho)\big)](1 -\beta).
\end{equation}
The lower bound \eqref{eq:4.57a} for RWRE now follows from \eqref{eq:4.60}, \eqref{eq:4.61} and Proposition \ref{prop:branching}.
Let us now turn to the proof of the lower bound \eqref{eq:4.57a} for ORRW. We follow the same strategy as above, and we see that \eqref{eq:4.60} and \eqref{eq:4.60a} hold. It remains to compute the probability of the event $C$:
$$
\mathbf P(C)=\Big(\frac b{b+\delta}\Big)^2\, \Big(\frac \delta{b+\delta}\Big)^{a/2-1}\,\Big(\frac \delta{b-1+2\delta}\Big)^{a/2-1}.
$$
By proceeding as in \eqref{eq:4.61} and above, and with the help of Proposition \ref{lemma:weights}, the proof of \eqref{eq:4.57a} is completed.
\qed \end{proofsect}
%\begin{remark} 
%The condition $\mathbb E[A^{-8}]<\infty$ can be relaxed to $\mathbb E[A^{-2-\varepsilon}]<\infty$ for some $\varepsilon>0$. This follows if we apply H\"older's inequality instead of Cauchy-Schwarz's inequality in \eqref{eq:3.7} and \eqref{eq:4.2}.
%\end{remark}

%%-----------------------------------------------------------------------------------------------------------------

\section{Appendix}
\small
%Let $\nu_{k,i},i=1,\ldots,b^k,$ be an enumeration of the vertices at level $k$, and  

\begin{Lemma}
\label{lemma:A.2}
Let $M(n,q)$ denote the $n$-th moment of a geometric random variable with parameter $q$. Then for $n \ge 1$, 
$M(n,q) \le c_n\, q^{-n}$, for some  constant $c_n$ that only depends on $n$.
\end{Lemma}
%For instance one can choose $c_1=1,\,c_2=2$.
\begin{proofsect}{Proof} 
We define  $g(q,n) \df \sum_{k=1}^{\infty} k^{n} (1-q)^{k-1}$, and notice that
$M_{n}^{\ssup q} = \sum_{k=1}^{\infty} k^{n} q (1-q)^{k-1} = q g(q,n)$. 
Since $0<q<1$, it is enough to show that there are coefficients $a^{\ssup n}_{\cdot}$ such that
\begin{equation}
\label{eq:7.1}
g(q,n) = \frac{\sum_{i=1}^{n} a^{\ssup n}_{i} q^{n-i}}{q^{n+1}}=\sum_{i=1}^{n} a^{\ssup n}_{i} q^{-i-1}. %\text{ with }a^{\ssup n}_{n} \le n!,\,\,a^{\ssup n}_{i} ???? \text{ if }  2 \le i \le n-1.
\end{equation} 
We prove  (\ref{eq:7.1}) by induction. As $g(q,1) = 1/q^{2}$, \eqref{eq:7.1} holds for $n=1$. Suppose now \eqref{eq:7.1} holds for $n-1$. We have 
\begin{equation}
\label{eq:7.2}
\begin{aligned}
g(q,n) - g(q,n-1) &= %\sum_{k=1}^{\infty} (k^{n} - k^{n-1})(1-q)^{k-1} = 
\sum_{k=1}^{\infty} k^{n-1} (k - 1)(1-q)^{k-1} \\
&= (1-q) \frac {\d}{\d (1-q)}  \sum_{k=1}^{\infty} k^{n-1} (1-q)^{k-1} = (1-q)  \frac {\d}{\d (1-q)}  g(q,n-1),
%&= (1-q)\frac {\d}{\d (1-q)}  g(q,n-1) 
%= \sum_{i=1}^{n}i a_{i}^{\ssup {n-1}} q^{-i - 1} -  \sum_{i=1}^{n}i a_{i}^{\ssup {n-1}} q^{-i}
\end{aligned}
\end{equation}
where  $ \frac{\d}{\d x}$  denotes the derivative with respect $x$. By the induction hypothesis,
\begin{equation}
\label{eq:7.3}
\frac {\d}{\d (1-q)} g(q,n-1)= \sum_{i=1}^{n-1}(i+1) a_{i}^{\ssup {n-1}} q^{-i - 2},
\end{equation}
and hence, using \eqref{eq:7.1} to \eqref{eq:7.3},
\begin{equation}
\label{eq:7.4}
\begin{aligned}
g(q,n) %&= \sum_{i=1}^{n-1}a^{\ssup {n-1}}_{i} q^{-i-1} + (1-q) \sum_{i=1}^{n-1}( i+1) a_{i}^{\ssup {n-1}} q^{-i - 2}\\
       &= n\, a_{n-1}^{(n-1)}q^{-n-1}+ \sum_{i=2}^{n-1} i(a_{i-1}^{(n-1)}-a_i^{(n-1)})q^{-i-1}-a_1^{(n-1)}q^{-2}.
%%g(q,n-1) + \sum_{i=1}^{n}i a_{i}^{\ssup {n-1}} q^{-i - 1}  -  \sum_{i=1}^{n}i a_{i}^{\ssup {n-1}} q^{-i}
%%=  \sum_{i=1}^{n-1} a_{i}^{\ssup {n-1}}\big(i(1-q) +q) \big) q^{-i - 1} 
\end{aligned}
\end{equation}
%Hence 
%$$
%a^{(n)}_n=n\, a_{n-1}^{(n-1)}, \text{ and for }2 \le i \le n-1,\,\,\,a^{(n)}_i=i(a_{i-1}^{(n-1)}-a_i^{(n-1)}),\,\,\,a^{(n)}_1=-a_1^{(n-1)}.
%$$
%Since $a^{(1)}_1=1$, 
This shows (\ref{eq:7.1}), and the proof is finished.
%moreover we see that $a^{(n)}_n=n!$ and $a^{(n)}_1=(-1)^{n-1}$....
%Now we use the induction step, i.e.  $a_{i}^{\ssup {n-1}} \le 0$, for all $i \le n-1$ and $a_{n}^{\ssup{n-1}} \le (n-1)!$. As $a_{i}^{\ssup n} =  a_{i-1}^{\ssup {n-1}}\big(i(1-q) +q) \big) $, we get that $a_{i}^{\ssup n} \le 0$ for $i \le n$ and 
%$$a_{n}^{\ssup n}  = a_{n-1}^{\ssup {n-1}} \big(n(1-q) +q) \big) \le a_{n-1}^{\ssup {n-1}}  n \le n!.$$
%To end the proof of Lemma~~\ref{lemma:A.2}, notice that
%$$ M_{n}^{\ssup q} = \sum_{k=1}^{\infty} k^{n} q (1-q)^{k-1} = q g(q,n) \le q \frac{n!}{q^{n+1}} =  \frac{n!}{q^{n}}.$$
\qed \end{proofsect}

\noindent{\bf Aknowledgement} A.~C.~was supported by the DFG-Forschergruppe 718  ÔAnalysis and stochastics in complex physical systemsÕ, and by the Italian PRIN 2007 grant 2007TKLTSR "Computational
 markets design and agent-based models of trading behavior". T.~S.~was supported by a postdoctoral research grant from the Max Planck Institute ``Mathematics in the Sciences''.

\end{document}